# RESTRUCTURING LOGIC

## Abstract

The outline of a programme for restructuring mathematical logic. We explain what we mean by 'restructuring' and carry out exemplary parts of the programme.

# Preface

A restructuring programme operates in three steps: First the existing materials of a discipline are liquefied, then passed through a filter, and finally allowed to settle and refreeze.

In the course of restructuring no substantive results are lost. (Some results will change their shape or their function, but none will be lost.) All "losses" – the dirt removed during filtering – lack mathematical substance. By shedding residues of idle metaphysics and obscure scholasticism we gain clarity and rigour.

At least half the work under the programme therefore consists of recycling and reframing existing results. This will substantially increase transparency in the areas processed. It will lead to new results only at a second stage, when the programme has opened up fresh vistas. Trickle down effects into other areas of mathematics are to be expected, but the main impact of the programme will be felt in the so-called foundations – in areas such as model theory, proof theory, recursion theory, and set theory.

Any restructuring programme is based on an abstract concept which acts as a solvent to liquefy the status quo. The quasi-philosophical position implicit in the status quo is opposed with an alternative position. The challenge inherent in the programme is to reconstruct the discipline on the basis of this alternative position.

In the present manuscript I will not be arguing philosophically for the positions at the heart of our restructuring programme. The philosophical arguments can take place elsewhere. Here I will simply state and, where necessary, explain the underlying philosophy. Instead of arguing for it philosophically we will let the fruitfulness of its mathematical consequences, and the clarity of mathematical presentation it allows, speak for itself.

The hegemonic if not entirely unchallenged position in mathematics is still a philosophy of Platonist realism about mathematical objects and a correspondence account of mathematical truth. This we reject.

Our alternative theory of truth falls largely outside the scope of this manuscript, although a few brief remarks will be made in the section on restructuring model theory. The alternative ontology is one of strict finitism – a conviction that all mathematical terms should ultimately be explicable in terms of finite operations on finite strings.

The idea of a mathematical reality is meaningful only if this reality is understood to be a projection of accepted theories. Formalism can be revived as the natural philosophy of mathematics provided it is embedded in an account of truth free from notions of correspondence. And provided, of course, the remaining Platonist habits that caused the first downfall of formalism are shed.

Certainty in mathematics is not certainty of correspondence. Certainty in mathematics means the certainty that a given formal sentence follows from a given set of formal axioms. These facts may be certain for all conceivable purposes, but they are not very interesting kinds of facts.

In mathematics the flow of justification is at least as much from consequences to axioms as from axioms to consequences. Axioms help organise consequences, they are usually designed after a batch of consequences is already known.

The image of mathematical logic as a foundation is therefore wrong. There is nothing particularly solid, or basic, or certain about the foundations of mathematics. On the contrary, results in certain areas of mathematical logic are necessarily aeronautical – thin, abstract, ethereal.

# Overview of Contents

The introductory section shows how model theory can be liberated from extraformal accretions.

Following that, because all serious objections to the programme arise from diagonal arguments, several versions of this type of argument are examined.

First to isolate the postulates and incompletely stated assumptions involved in each version.

Second to argue that greater clarity, simplicity and power can be achieved by making changes to these postulates: In some cases by rejecting traditional postulates, in other cases by admitting alternative postulates in competition to the existing ones. (On the table of contents above, sections where rejection is mandatory for the programme to succeed are marked by 'M', sections where it is optional by 'O'.)

# Restructuring Model Theory

What is the purpose of a model ? To assign one of two values to every sentence of a theory in a way that respects logical structure. The best understood and most transparent way of doing this is a decidable theory. Since we expect our models to behave like theories, we might just as well insist that they ought to be theories.

There are objections, of course. Model theory in its current form is plagued by various extraneous commitments. First there are commitments to diagonal arguments that make the idea of the model as a theory seem problematic. Depending on how it is viewed, the idea of an uncountable model that arises out of Cantor's argument can appear hard to reconcile with countable typography. And then there is the irredeemable insufficiency of theories alleged by Gödel's diagonal incompleteness argument. It appears to imply that if models had to be decidable theories, only weak, uninteresting theories could have models.

These non-philosophical commitments are dealt with in other sections, here we will only very sketchily address the more philosophical ones before proceeding to a new definition of models.

Model theory cannot help entering into philosophical commitments – because of its declared objective to give a definition of truth appropriate for mathematics it is inevitably the dressing up of a philosophical theory of truth in mathematical symbols. Whatever logicians may like to claim, it is impossible for model theory to be philosophically neutral.

When theories fracture – as they would appear to do under the current interpretation of Gödel's theorem – one logical response is to seek a whole outside of theories. (The other typical response is a defeatist retreat to some form of relativism.) Gödel himself was driven, arguably in part by his own theorem, to a philosophical position that stresses the reality and independence of mathematical objects divorced from theories. He became known as the arch-Platonist. This is a position that chimes well with model theory in its current form. The spirit of current model theory is one that breathes a commitment to a correspondence theory of truth.

Hence a choice of two philosophical commitments for model theory and thereby for mathematics as a whole:

> **The traditional favourite:** Truth consists in correctly representing an external state of affairs.

For this position, theories are secondary, in the extreme only faint reflections of quasi-Platonic perfect forms.

Due to constraints on space the alternative position can only be floated in the form of a slogan.

> **The alternative:** Language is a system of presentation, not representation.

Since formalised systems are by far the most transparent and trustworthy way of presenting mathematics, this leads to neo-formalism, or the thesis that nothing outside of theories is admissable: Every mathematical argument ought to be explicable in terms of either derivations within or relations between theories.

To qualify as true a theory must not only be well-engineered – consistent and complete –, but also attractive and perhaps even useful. There is a beyond to the formal but this beyond has little to do with correspondence to an external reality. One element of the beyond is perfectability, a constant exhortation to make existing theories better; the other is asking for the point of it all this activity.

The possibilities of devising formal games are endless. The question that truth asks is one that goes outside of model theory: Are they *worth* playing ? Do they make *good* mathematics ? Whether a theory is true in the colloquial sense amounts to a toned-down value judgement that is more moral or aesthetic than factual.

Apart from the objections, there are also clear and present benefits to be derived from identifying models with decidable theories. Model theory in its current form degrades formality: Truth for theories is defined in such a way that the models theories are referred to are less formal than the theory.

Current model theory takes axiomatised theories which represent everything that is desirable in mathematics – clear, clean, and crisp – and turns them over to a murky so-called model or interpretation, a semi-formalism at best, where notation is an unstandardised mishmash with usually a large helping of naïve set theory. Applying the alternative standard of truth, current models must seem considerably less real than theories. The operation of interpretation the way it is conducted today constitutes a crime against clarity.

Completeness of a theory is a much clearer demand than the law of excluded middle currently assumed for models; consistency of a theory much clearer than the existence of a model, which whatever the real existence of models was supposed to mean surely meant to imply consistency. The issue is not whether consistency and completeness assumptions are made in model theory, but how transparently it is done.

The semantics of current model theory only give the illusion of reaching "real mathematical objects". What is reached in reality are only unformalised notions about such objects, i.e. a less explicit, inchoate type of theory. The issue is not one of reaching through to a supposed reality rather than staying with symbolic presentations, but of how well articulated the symbolic presentations are.

Current semantics thus fail to add value. In reality all they do is to provide is an excuse for evading the challenge of demonstrating the coherence (consistency) and definiteness (completeness) of unformalised notions.

When current model theory allows models to be defined by pointing to unformalised notions, it effectively tries to deflect the challenge of making them clear.

Relating a formal system to the world is less straightforward than the image of mirroring wants to make us believe. One should not think of filling an empty formalism with prefabricated objects, but rather of attempting to draw out a formalism from an unformalised account of the world. The formal assimilates the less formal. Finding a formal system "exemplified" or "instantiated" is a complicated process of trying to distil a less ambiguous symbolic presentation out of a more ambiguous one in such a way that the end product still "reasonably resembles" the original input. (Note: The vagueness of phrases like 'reasonably resembles' is not accidental, but to some extent inevitable and therefore appropriate.)

Only by way of a demonstration that there are plausible formal theories to be extracted from it can an unformalised notion prove its worth. Where no demonstration is forthcoming we have a right, even the duty, to remain cynical.

By making it clear that model theory has no concern for anything outside of the formal, that its subject matter are the provability invariant transformations of theories, we allow model theory to join the mathematical mainstream for perhaps the first time. Nothing could be more in accord with the mathematical mainstream than studying transformations and their invariances – it would hardly be an exaggeration to say that every branch of mathematics could be defined as the study of this-or-that-invariant transformations –, nothing more alien to it than the models of current model theory, which can look like ideas about "real" objects that have come only part of the way towards refinement into a theory.

It is time for model theory to shed its residual commitments to correspondence thinking and the attending substandard formalisms. Much better to use definitions that prevent a loss of formal clarity in moving from the theory to the model.

> **Definition:** A modelling transformation is a function m from the strings of a theory Th to the strings of a theory M such that
>
> 1.  m respects the roles of strings, mapping
>
>> terms of Th to terms of M,
>> predicates of Th to predicates of M,
>> sentences of Th to sentences of M.
>
> 2.  m commutes with logical connectors, quantification, and substitution, i.e.
>
>> for all sentences S,           M |- m(¬S) ↔ ¬m(S),
>> for all sentences S, T         M |- m(S ∨ T) ↔ m(S) ∨ m(T),
>> for all predicates P,          M |- m(∃x P(x)) ↔ ∃x m(P(x)), x free,
>> for all predicates P, terms t  M |- m(sub(P,t)) ↔ sub(m(P),m(t)),
>
> 3.  m preserves provability, if Th |- S then M |- m(S).

For systems based on functions and equations rather than predicates and sentences an equivalent definition can be given, essentially only by introducing the clause

> for all functions f, terms t    M |- m(sub(f,t)) ) ≡ sub(m(f),m(t)),
> with ≡ the image of = under m

> **Proposition:** The image of = under a modelling transformation is an equivalence relation.
> Proof: By elementary logic.

The new definition relates to a more familiar one by

> **Proposition:** For m to be a modelling transformation it is sufficient that
> 1.  Every atomic relation from Th maps to a relation in M of the same arity (not necessarily atomic)
> 2.  Every primitive term from Th maps to a closed term (not necessarily primitive)
> 3.  Equality is mapped to equality
> 4.  Non-atomic predicates are mapped recursively according to
> m(¬S) = ¬m(S),
> m(S ∨ T) = m(S) ∨ m(T),
> m(∃x P(x)) = ∃x m(P(x)),
> m(sub(P,t)) = sub(m(P),m(t)).

5. Non-primitive terms are mapped recursively according to
   m(sub(f,t)) ) = sub(m(f),m(t)).
6. m(A) is provable for all axioms A of Th
   Proof: By induction on the complexity of formulae.

Given target values for the primitives all strings of the theory can be parsed and rebuilt from bottom up in the model.

**Definition:** If Th maps into a consistent, complete theory M, then M is said to be a model for Th.

**Definition:** A sentence S of Th is said to be true in a model M if M |- m(S), the image of S under the modelling transformation.

Exactly how these definitions are phrased is in many details a matter of convenience. Some people might, for example, prefer to refer to what is defined above as a para-model, and reserve the term 'model' for a stricter concept that includes a condition of maximal satisfiability. Our concern here is less with establishing any particular definition than with outlining minimal conditions for a satisfactory definition of a model:

- The modelling relation – the relation between theory and model – should be single-typed. The model must be the same kind of thing, ontologically, as the theory. Domain and codomain of the modelling relation should share the same range: strings of a language over a certain logic.

- All relations in model theory must ultimately be reducible to finite operations on strings. Even when we talk of structures as equivalence classes of theories no non-formal realities should be implied: A concrete theory is taken as a token, and for another theory to be in its equivalence class is to be mappable to this theory. Mappability being a finite string relation.

- No shifting of orders. Predicate expressions are to be mapped to (second order) predicate expressions, not to (first order) object expressions. Even when the model is in set theory the image of a predicate expression must be a set-theoretic predicate. (This set-theoretic predicate may in a second step be extensionalised – hence converted to a first order object – but secondary conversion does not change the fact that the modelling transformation maps to the predicate, not the extension.)

- No precommitting to set theory. Model theory should start with a definition that refers only to mappings from strings to strings. These mappings could turn out to have equivalents in set theory, but these parallels should be considered to be purely coincidental. There should be no use or mention of (infinitistic) set theory in the basic definition of a model.

- No hint of shifting to a different ontological plane, exchanging strings for objects. The modelling relation should be defined unambiguously as a transformation among equals, a way of circling around, without any pretence of exchanging formal for non-formal objects. This is an approach naturally compatible with category theory. The definition should not, however, precommit to category theory as a foundational framework, either.

- No use of a concept of truth or appeal to correspondence relations. While it is possible to interpret much current model theory syntactically – and many logicians in their folkways in effect already accept highly syntactic definitions – there is nothing in the standard textbook definition of a model to preclude interpreting the term 'model' as a structure in the mind of God. This ambiguity is embarrassing, and has to end.

All of these conditions clearly are satisfied by our new definition, and just as clearly not satisfied by the traditional definition going back to Tarski. To prove that the difference is real we give a first example of a word game that would not have been possible under the traditional definition:

**Proposition:** A complete and consistent theory is its own model (under the identical mapping).
Proof: String identity implies equivalence.

Restated in slightly different terms, the proposition says that every theory has at least one automorphism. In fact, homomorphisms everywhere in mathematics can be viewed as modelling transformations. The usual results about homomorphisms apply. In particular, the automorphisms of a theory form a group.

**Definition:** A sentence is said to be valid if it is provable in the axiomless predicate calculus.

**Proposition:** A valid sentence is true in all models.
Proof: A valid sentence is provable in any theory based on the predicate calculus, including all models.

Note how validity – and as we shall see next, implication – obtains only relative to a specified range of models. Results would be different if, for example, second order models were allowed.

**Definition:** A sentence S is said to imply a sentence T if in all models M, M |- m(S) then M |- m(T).

**Proposition:** If Th |- S → T, then S implies T.
Proof: By modus ponens in M.

The task now is to reduce inessential typographical variation.

**Definition:** Two modelling transformations $m_1$, $m_2$ are elementarily equivalent if they declare the same sentences true, i.e. $M_1$ |- $m_1$(S) iff $M_2$ |- $m_2$(S).

**Proposition:** The models of a decidable theory are elementarily equivalent.
Proof: By definition.

At the other extreme, imagine a set of sentences without axioms. Every elementarily non-equivalent valuation though a model corresponds to a different choice of axioms. How many different ways are there of turning a set of sentences into a decidable theory ?

**Definition:** The predicates of a decidable theory Th are the equivalence classes of its predicate strings under the relation: P≡Q iff Th |- ∀x P(x) ↔ Q(x).

**Proposition:** Two modelling transformations are elementarily equivalent if they agree on the equivalence of predicate strings.
Proof: Suppose the transformations disagree on at least one sentence. This sentence is equivalent to a sentence of the form Qx P(x). Then P must be equivalent to either a constant true or a constant false predicate under one transformation and non-equivalent under the other.

A theory is thus determined by the way in which it sorts its predicate strings into equivalence classes. (Not for nothing is the logic called the predicate calculus.) When transforming theories, this allows us to concentrate on how the predicates are mapped.

**Definition:** A modelling transformation is said to be surjective (bijective) if it is surjective (bijective) as a mapping of predicates.

**Proposition:** The inverse of a modelling transformation m: Th → M is a modelling transformation in its own right when Th is decidable and m surjective.

**Corollary:** A surjective modelling transformation between decidable theories is bijective.

**Definition:** A theory isomorphism is a bijective modelling transformation.

**Definition:** A structure is an equivalence class of theories under the relation of theory isomorphism.

The best referent for the phrase 'the natural numbers' is not a set, or any collection of objects, but a structure in this sense – number theory up to typographical variation.

**Definition:** A theory is categorical if its models are isomorphic to each other.

Strictly speaking, a theory is categorical only for a logic, e.g. the finite order predicate calculus. (Distinguishing between finite orders makes little sense as set theory cuts across orders.) By claiming that a theory is categorical *tout court* one would be predicting a limit for human ingenuity – always a hazardous undertaking. It is at least possible that in future more powerful formalisms could be invented that will be able to draw finer distinctions than any currently known.

**Definition:** The objects of a decidable theory Th are its atomic predicates, i.e. the equivalence classes of predicate strings such that Th |- ∃!x P(x).

The objects of a theory are implied by how discerning its predicates are. When a language lacks the means to draw a particular distinction, objects constituted by this distinction do not exist for the language.

The definition is a sensible one because as soon as a theory is able to identify its objects uniquely whatever is true for all 'objects' in the sense given will be true universally:

**Definition:** A theory is said to be atomic if whenever Th |- ∃x P(x) there exists a predicate Q such that Th |- ∃!x P(x) ∧ Q(x).

**Proposition:** For an atomic theory it is provable (in sufficiently powerful models) that

$$\forall P \; ( \; \forall Q \; \exists!x \; Q(x) \rightarrow ( \; \forall x \; Q(x) \rightarrow P(x) \; ) \; ) \; \rightarrow \; ( \; \forall x \; P(x) \; )$$

**Proposition:** Categorical theories are atomic.

The propositions suggest that constants can be eliminated in favour of predicates, but not the reverse. Traditional model theory takes a different view. It wants to step out of transformations and get to objects, which are believed to pre-exist, not constituted by the non-identities theories are able to discern.

Hence traditional model theory is more or less committed to a project of extensionalisation, or the elimination of predicates in favour of objects. Set theory lends itself to this project, as only set theory could easily facilitate a conversion of second order predicates into (sets of) first order objects. Models so far have essentially been models in set theory, be it models in an inarticulate version of set theory with urelements and tilt towards the naïve.

While it is possible to find models in set theory under our paradigmatic definition of a modelling transformation, the more common extensionalising transformations are enabled only by this generalisation:

**Definition:** A generalised modelling transformation between formal systems Th and M is a mapping from the strings of Th to the strings of M that preserves
>> the roles of strings,
>> derivability,
>> and commutes with the construction principles,

where the construction principles are string mappings in Th and M such that
>> substitution: terms × open sentences → sentences
>> disjunction: open sentences → open sentences
>> negation: open sentences → open sentences
>> quantification: open sentences → sentences

The tight mapping from predicate strings of Th to predicates of the same order in M remains available under this definition. In addition, however, looser modelling transformations are introduced that map predicates to (expressions for) subsets of the universe, or even (expressions for) indicator functions. The substitution function takes on a different meaning as well.

$$m(P(a)) = m(sub_{Th}(P,a)) = sub_{Set}(m(P),m(a)) = sub_{Set}(p,a) = a \in p$$

The traditional idea of a model is perhaps best understood as representing a special case of such a (loose) modelling transformation. It transports a theory to an image in set theory.

**Definition:** The Tarskian image of a theory in set theory is an existential statement produced as follows:
1. Without loss of generality assume that Th is finitely axiomatisable in either first or second order logic.
2. Eliminate constants and functions in favour of relations.
3. All axioms can then be condensed into a single axiom for which we write as F($R_1$, ...., $R_n$), with relations $R_i$
4. The translate is then the statement
   "∃u,$r_i$ such that u is a set (the universe), and the $r_i$ are relations over u (subsets of the appropriate Cartesian products) such that F*($r_1$, ...., $r_n$)", with x∈$r_i$ replacing $R_i$(x) in F*.

Due to its capacity to receive an image of any finitely axiomatisable theory set theory should perhaps be regarded more as a scheme for defining theories than an individual theory.

**Definition:** A theory is said to have an extensional model if its Tarskian image is provable in set theory.

One of the basic assumptions traditional model theory has made is that the existence of an extensional model (in whichever tacitly assumed version of set theory) should be equivalent to the consistency of the theory. One half of this is true.

**Proposition:** A theory is consistent if has an extensional model.
Proof: The image was defined to ensure that any inconsistency in the theory carries over into set theory.

The converse is false. In actual fact, consistency is weaker than the existence of an extensional model. There are theories that are consistent, yet non-extensionalisable. Set theory is one example.

The idea of extensionalisation in model theory is therefore only partially salvageable, and only if interpreted strictly as a syntactic modelling transformation from a theory into set theory, with no atomic individuals apart from sets.

So while the model in set theory might occasionally be useful, it is only an optional extra. The primary model for an algebraic theory, say, will always be an algebraic theory.

In model theory, restructuring does not directly provide new results. It greatly improves presentation, and could suggest fresh lines of inquiry that might lead to new results in the future, but does not profoundly change current practice.

The next section is going to be different. There existing results will undergo significant changes.

# Restructuring Gödel's Incompleteness Argument

## I.    Gödel's First Incompleteness Theorem
### or
### The Alleged Incompleteness of First Order Number Theory

**N.B.** Adopting the changes advocated in this section is a *conditio sine qua non* of the programme. The traditional interpretation of Gödel's diagonal argument, if upheld, would trivialise the new definition of a model.

Gödel's argument ingeniously shows that some theory has a model in number theory. This section is devoted to answering the question: Which theory ?

The pre-image for Gödel's mapping of strings into number theory must be a theory in which the strings of number theory figure as closed terms. String predicates – e.g. the predicate of being a well-formed sentence – must first exist, unencoded, in a string theory before they can be projected, under a code, into number theory. The theory of number-theoretic strings is where string predicates live, in number theory they are only guests.

> **Definition:** A theory M is said to be meta to a theory Th if all strings over the alphabet Σ of Th are well-formed as closed terms in M.

To assume the existence of a set of strings is to assume at the very least the existence of a predicate in a decidable theory of strings. It would be irresponsible to trust that talk of sets of strings, especially infinite sets, is in any way coherent unless and until such a theory is presented.

So before Gödel's encoding argument can even begin, the following assumption has to be made:

> **Gödel's Assumption:** There exists a decidable theory META that is meta to number theory and contains a predicate Prov() that correctly expresses provability for number theory, i.e. for all sentences S of number theory $\mathbb{N}$,
>
> > META proves/does not prove Prov(S) iff $\mathbb{N}$ proves/does not prove S.

Note that the strings of $\mathbb{N}$ appear in META nakedly, untouched by any encoding.

> **Lemma:** Gödel's assumption can be analysed into
>
> > **Gödel's First Assumption:** Concatenation is formalisable. In other words, there exists a theory CONCAT containing a binary function c(,) that correctly expresses concatenation over a given finite alphabet Σ, so that in particular, for all primitive terms a,b,c∈Σ*,
> >
> > > CONCAT |- c(a,b) = c iff c is the concatenation of a and b.
> >
> > **Gödel's Second Assumption:** Recursive definition is sound, i.e. extensions of CONCAT by predicates that can be defined recursively in terms of concatenation are consistent if only CONCAT is.
> >
> > In particular, the extension of CONCAT by the following seven relations, all of which are definable in terms of concatenation, is consistent:
> >
> > > VARIABLE, OPEN_TERM, CLOSED_TERM, SUBSTITUTION, PREDICATE, SENTENCE, PROOF.
> >
> > **Gödel's Third Assumption:** The theory consisting of CONCAT and the seven relations necessary to define PROOF has a complete extension.
>
> Proof: See Section Appendix

The first assumption is infuriatingly imprecise. 'Correctly expresses' is just a pretentious way of saying 'behaves as one would naïvely expect concatenation to behave'. But the task we face is impossible: How does one clearly state an incoherent intention ? How do we even know that there is something there to be expressed ?

In order to get a firmer grip let's try

> **Definition:** CONCAT$_n$ is a theory that correctly expresses concatenation for strings up to length n over a given alphabet.

Writing out examples of CONCAT$_n$ is straightforward if tedious. To illustrate, we give a simplified version of CONCAT$_2$, a decidable theory consisting of three primitive terms **0**,**'**,**0'**, a single predicate letter C, and the following axioms:

| | |
|---|---|
| $\forall x \; x = \mathbf{0} \lor x = \mathbf{'} \lor x = \mathbf{0'}$ | There is an explicit, pre-declared finite domain. This in itself is nothing unusual, similar limitations are implicit in the workings of any real world computer. |
| $\forall x,y \; x \neq \mathbf{0} \lor y \neq \mathbf{'} \rightarrow \neg\exists z \; C(x,y,z)$ | Concatenation in CONCAT$_n$ is only a partial function, not defined beyond a given upper limit (again analogous to real world computing). Limit setting could be made more sophisticated by appealing to string length rather than explicit listings. |
| $C(\mathbf{0'},\mathbf{0'})$ | Within the strictures of the predicate calculus the only way of giving correct extensional values is to list all instances individually. |

It seems reasonable to grant the intelligibility of the idea of constructing theories CONCAT$_n$ on this pattern even for n larger than the physical capacity of any human or computer.

With the sequence CONCAT$_n$ in mind we are now able to clarify

> **Gödel's First Assumption (refined):** The sequence of languages CONCAT$_n$ has a limit, i.e. there exists a decidable theory CONCAT such that the theories CONCAT$_n$ have nested models in CONCAT.

> where

> **Definition:** A sequence of languages L$_n$ is said to have a nested model in L if for all n,

>> L$_{n+1}$ is a model for L$_n$, and

>> sentences true in L$_n$ are true in L for strings up to length n.

> **Corollary:** L is a model for all the L$_n$.

Once these assumptions are in place it becomes possible to connect META with number theory by way of two lemmas.

> **Embedding Lemma for Meta Theories:** META$_{Th}$ contains a model for any decidable theory Th.
> Proof: Let predicates P(x) map to Prov(sub(P,x)), and terms identically. The mapping describes a model provided Th is decidable. (It needs to be shown that META |- Prov(sub(¬P,x)) ↔ ¬Prov(sub(P,x)). This is true if Th is decidable, and since Prov() expresses, provable.)

> **Gödel's Lemma:** ℕ, if decidable, contains a model for any consistent subset METAPART of META.
> Proof: This is a restatement of Gödel's groundwork for the First Incompleteness Theorem, accepting the substance with only minor changes: Instead of mapping (elements of) sets of strings into number theory, predicate strings from META are mapped to predicate strings of ℕ. The consistency assumption, although not stated, is implicit in Gödel's argument.

> **Definition:** A proper subtheory of a theory Th that is a model for Th is called an inner model.

> **Corollary:** If ℕ is decidable, META contains an inner model.

Let g map META into $\mathbb{N}$ [†], let f map $\mathbb{N}$ into META, and for sentences S put S* = f(g(S)).

We now have all the ingredients for proving that the predicates of META under f o g have a fixed point. This follows generically from, for instance, category theory. We will nonetheless give the proofs in long form.

**Round Trip Lemma**: META |- S ↔ S*.

Proof: META |- S implies $\mathbb{N}$ |- g(S), which in turn implies META |- f(g(S)). Therefore META |- S*.

For the converse suppose META |- ¬S. Then $\mathbb{N}$ |- g(¬S). So because of modelling, $\mathbb{N}$ |- ¬g(S). Then META |- ¬S*.

The cyclical image S* of S in META must be if not equal then at least equivalent to S.

**Definition:** The predicate True is a truth predicate for a meta theory META$_{Th}$

if for all sentences S, META |- True(g(S)) ↔ S,

where g is a modelling transformation from META to Th.

**Gödel's Theorem (according to Gödel):** If $\mathbb{N}$ is decidable, META contains a truth predicate.

Proof:   We show that Prov() is a truth predicate.

META |- Prov(g(S)) ⇒ META |- S.

As Prov() expresses provability for $\mathbb{N}$, META |- Prov(g(S)) implies $\mathbb{N}$ |- g(S).

By embedding, $\mathbb{N}$ |- g(S) implies META |- f(g(S)) = S*.

Hence META |- S (Round tripping).

META |- Prov(g(S)) ⇐ META |- S

META |- S implies $\mathbb{N}$ |- g(S).

As Prov() expresses provability for $\mathbb{N}$, META |- Prov(g(S))

The plot thickens with

**Tarski's Theorem:** No decidable theory can contain a truth predicate.

Proof: By applying to itself the predicate Fssb(P) := ¬True(g(sub(P,g(P)))), (F for false, ssb for "self-substitution").

**Corollary:** Either

Gödel's assumption (= META exists and is decidable)

or that

Number theory is decidable

is false.

Proof: Taken together the two would prove the existence of a decidable theory META with a truth predicate, contradicting Tarski's Theorem.

At first sight there appears to be an honest choice. The question seems to be on which of the initial assumptions should we pin the blame for a contradiction that emerges at the end of a long derivation. The string theory

---

META is a prima facie favourite for carrying the blame – it is an odd, untested and hugely ambitious theory. This conclusion becomes inevitable once we note that essentially the same argument can be made without involving number theory at all, simply by taking a second copy of META in place of number theory. (It goes without saying that the argument can *not* be made with two copies of number theory.)

**Gödel's Theorem (2nd approximation):** $META_{Meta}$ contains a truth predicate.
Proof: Let $META_{Meta}$ be a variation of META, the meta theory for a concatenation-based meta theory.
Take two copies of $META_{Meta}$, labelled $META_1$ and $META_2$. Each theory is, without any circularity, a meta theory for the other.
One way still takes the Embedding Lemma, the way back is much easier this time: the identical mapping. Round Trip then holds true, so the Theorem follows.

**Corollary:** $META_{Th}$ exists for no theory that is free of limitations on string length.
Proof: $META_{Meta}$ is in no way special. Any variation of META for different underlying theories consists of a version of CONCAT and extending predicates. The differing versions of CONCAT are models of each other by a simple permutation of the alphabet; the meta predicates for one formal theory are definable by the same means as the predicates for any other.

Gödel's incompleteness argument, though revealing about the string theory, tells us nothing about number theory. The contradiction used to refute the decidability of number theory is frankly imported from META. Given Gödel's original conclusion was the result of projecting the failures of a string language into number theory, no reason remains for believing in the undecidability of number theory. Although absence of counterevidence does not constitute evidence, it is only natural to revert to decidability as the default assumption.

**Gödel's Theorem (3rd approximation):** There is no universal string language, i.e. no decidable theory that contained CONCAT and allowed recursive definition.

META is a myth. Ideas of such a theory people may have won't fly. There is nothing solid under the vapour.

So the real choice is given by

**Corollary:** One of the following three must be false

"CONCAT exists" (Gödel's First Assumption)

"Recursive definition is sound" (Gödel's Second Assumption)

"Any consistent theory has a decidable extension (and therefore a model)" (Zorn's Lemma)

Proof: Assume CONCAT existed. Then if recursive definition is sound, CONCAT extended by the predicates necessary to define provability is consistent. Then if any consistent theory can be extended to a decidable theory, META exists and is complete. Contradiction.

While some may be prepared to doubt the consistency of number theory, the admissability of recursion, or the validity of induction as a rule of inference – all of which involve Gödel's second assumption –, and others may want to doubt that all theories have complete (and by definition: axiomatisable) extensions, for the rest of us this is as good a *reductio ad absurdum* of CONCAT's existence as one can hope for.

**Gödel's Theorem (4th approximation):** $\neg\exists$ CONCAT.

The lesson of Gödel's Theorem: Concatenation, used naïvely, can be just as treacherous as the membership relation $\in$. So in a way Gödel's argument does for string theories what Russell's paradox does for set theory. For both theories, the naïve assumption that all predicates ought to be able to have extensions turns out to be inconsistent despite its overwhelming intuitive appeal.

**Corollary:** For any decidable theory Th there exists an n such that $CONCAT_n$ has no model in Th.
Proof: The theory would contain a model for CONCAT if $CONCAT_n$ embedded for all n.

The definition of meta theories now needs to be adjusted to the realities. All that can be hoped for is an incomplete theory scheme rather than the complete theory META would have been.

**Definition:** The theory PROOF THEORY is as META defined earlier, but with CONCAT replaced by the following axioms:

1) "Letters are primitive, i.e.

$\forall x\ (\ \text{LETTER}(x) \leftrightarrow \neg\exists u,v\ c(u,v) = x\ )$

with LETTER(x) defined as $x = a_1 \lor x = a_2 \lor ... \lor x = a_n$, for $a_i \in \Sigma$"

2) "Right- and Left-Concatenation are one-to-one, i.e.

$c(x,u) = c(x,v) \lor c(u,x) = c(v,x) \rightarrow u = v$

3) "The universe consists of strings arising through concatenation, i.e.

$\forall P\ (\ (\ \forall x\ \text{LETTER}(x) \rightarrow P(x)\ ) \land (\ \forall x,y\ P(x) \land P(y) \rightarrow P(c(x,y))\ )\ ) \rightarrow \forall x\ P(x)$

What, after all, is successor ? Just a function letter in a formal structure. So if we can define successor implicitly – as we must in order to reach infinity – then there is some hope that concatenation could be largely defined by the role it performs in relations leading up to provability.

Unlike number theory, however, PROOF THEORY contains relations that fail to live up to the meaning adumbrated in extensional approximations.

**Definition:** For a predicate P from PROOF THEORY let $P_n$ be the sequence of languages that correctly approximate P up to strings of increasing length n.

Evidently, the sequence $P_n$ for P is a generalisation of $CONCAT_n$ for concatenation.

**Definition:** A relation in PROOF THEORY is said to be extensionalisable[‡] if it can serve as the limit for its approximating sequence $P_n$.

**Gödel's Theorem (final):** There are non-extensionalisable predicates in PROOF THEORY.

Concatenation is not the only non-extensionalisable relation. For obvious reasons the same applies to the combination of provability and substitution. The universal string language in which all of these relations might have been extensionalisable can thus be thought of as the non-existent limit of a sequence of languages that grows out of bounds.

**Definition:** Let $PROOF\ THEORY_n$ be PROOF THEORY extended by $CONCAT_n$.

$PROOF\ THEORY_n$ can alternatively be arrived at by replacing CONCAT in META with $CONCAT_n$. If all predicates were extensionalisable, $PROOF\ THEORY_n$ would converge, turning an extension of PROOF THEORY into META.

So the General Embedding Lemma for meta theories fails. But for specific theories it can succeed, and specific predicates can be extensionalisable.

**Embedding Lemma for Number Theory:** $\mathbb{N}$ is isomorphic to a subtheory of PROOF THEORY.
Proof: Successor is modelled by concatenation with ', when the universe of strings is restricted to numerals of the form $0^{(n)}$. (For more details see Section Appendix.)

**Corollary:** The function c(x,**'**) is extensionalisable.
Proof: It is isomorphic to successor.

The lemma shows how arithmetic can be conducted with numerals. (To get from numerals to "numbers", a numeral like 0'' would have to be taken as a primitive term rather than a nested functional expression s(s(0)).)

---

[‡] This definition of 'extensionalisable' is equivalent to the set-theoretic meaning: Corresponding to every stage $P_n$ of the sequence we can define the finite set of all strings of length less than n that provably satisfy P(x). The infinite union, or limit, of these finite sets would then be the extension of P, the set of all strings such that P(x).

**Corollary:** Number theory is at most as powerful as PROOF THEORY.

with

**Definition:** A theory Th is said to be at most as powerful as a theory M if there exists a modelling transformation from Th to M.

For the ranking by power to be effective, M does not necessarily have to be decidable, hence a model, for Th. Modelling transformation can exist also between incomplete theories. Top of the ranking are the inconsistent theories, least powerful the empty theories (theories that prove not a single sentence).

**Corollary from Gödel's Theorem:** PROOF THEORY is strictly more powerful than number theory.
Proof: Its non-extensionalisable predicates behave unlike any predicates in number theory.

The greater power of PROOF THEORY is owed to what are in effect empty predicate names.

**Definition:** A predicate Q in a model M is said to represent a predicate P from a theory Th if Q = m(P) under the modelling transformation m: Th → M.

**Gödel's Lemma (revisited):** Every extensionalisable predicate in PROOF THEORY is represented by a predicate in number theory.
Proof: Gödel's original idea of mapping sets of strings into number theory can only be realised when the predicate strings that supposedly denote these sets are extensionalisable.

So meta theories are more powerful than the underlying theory but not for the reason given by the General Embedding Lemma, which had claimed that already a small portion of the meta theory – the sentences of the form Prov(sub(P,t)) – could express the underlying theory. The reason meta theories are more powerful is that they contain a compartment for empty predicate names, of which Prov(sub(x,y)) is one.

This then is Gödel's paradox: Concatenation-based string theories are meta theories for any reasonable language, including themselves. Diagonalisation proves that a fully extensionalisable meta theory must be more powerful than the underlying theory. So because meta to itself, a universal string language would have to be more powerful than itself. Contradiction, or the ordering by power breaking down.

# II.    Gödel's Second Incompleteness Theorem
## or
## The Confirmed Unprovability of Consistency

The absence of a universal meta language induces an asymmetry with respect to proving decidability.

The positive cannot be proved in an ordinary sense. It seems inescapable that the theory M *in* which we prove either consistency or completeness must contain a model of the theory Th *for* which we prove it, and so must be at least as powerful. This is because M can only meaningfully claim to be 'about' the underlying theory if the full complexity of Th carries over into the proof in M. With a benign twist of circularity, what is proved must be assumed. The assumption of decidability is never reduced, only transformed into equivalent forms – which can nonetheless be enlightening.

The negative is somewhat different. It is possible to *dis*prove decidability for particular versions of a formalism. But once the latest instance of inconsistency or incompleteness is remedied by a modification in the theory, openness reasserts itself.

So although not demonstrable, decidability is in a broad sense falsifiable. In other words, the claim of consistency for a theory constitutes a prediction and must be warranted by induction from experience. The sun will rise again tomorrow, and no one will find an inconsistency in number theory.

The expectation of reductive proofs of consistency has always been unreasonable. If by finitary one were to mean for instance a proof in a bounded quantifier theory then to expect that number theory could be proved consistent by finitary means is to expect number theory to have a model in a bounded quantifier theory, and so to expect that the use of unbounded quantifiers could be shown to have been spurious and inessential.

No means of proof can show that more powerful means are sound (This is what it must mean to be 'more powerful', according to any sensible definition). We can shift our ignorance about which sentences one theory will prove into ignorance about another. We can build transformations between theories. What we cannot do is to find foundations: The belief that decidability, if not already achieved, is achievable precedes all work in mathematics. In practical terms, the decidability of at least the oldest and best tested theory – number theory – is taken for granted.

This belief in being able to make theories work is arguably the one useful meaning of the otherwise misguided doctrine of mathematical Platonism, the idea of a consistent and definite mathematical reality existing prior to theory design.

## Section Appendix

**Proof of Lemma** – Gödel's First, Second and Third assumptions together imply for any theory Th the existence of a decidable theory META$_{Th}$ that contains a predicate Prov() expressing provability for Th.

> **Sublemma:** The relations VARIABLE, OPEN_TERM, CLOSED_TERM, SUBSTITUTION, PREDICATE, SENTENCE, PROOF are definable in terms of concatenation.
> Proof: Specifics in the definition of meta relations depend on the underlying theory. We will give definitions for one particular formulation of number theory. The variations, however, are inessential.

A typical example of a meta relation is NUMERAL() (as they are called in number theory, more generally one would speak of CLOSED_TERMs).

It is defined by the axioms

NUMERAL(**0**)                                           ("0 is a Numeral")

$\forall$x NUMERAL(x) $\rightarrow$ NUMERAL(c(x,'))        ("If X is a Numeral, then X' is a Numeral")

$\forall$P ( P(**0**) ) $\wedge$ ( $\forall$x NUMERAL(x) $\rightarrow$ P(c(x,')) ) ("Nothing else is a numeral")
   $\rightarrow$ ( $\forall$x NUMERAL(x) $\rightarrow$ P(x) )

The third axiom is necessary to exclude from NUMERAL everything that did not originate from its generation rule.

It is not hard to see how the predicates VARIABLE(x), OPEN_TERM(x), PREDICATE(x), SENTENCE(x) can be defined by a similar trinity of a base, a generation, and a closure clause. This leaves only two more relations – SUBSTITUTION(x,y)=z, which must be primitive recursive in concatenation if only the underlying theory is formal, and PROOF.

Intuitively, a proof is a string of sentences such that each sentence follows from its predecessor(s). From the fact that $\mathbb{N}$ is formal, we have that the relation FOLLOWS_IN_ONE_STEP(,) is primitive recursive and definable in terms of concatenation. We also assume that there exist a predicate codifying AXIOM(), which need not be more than an explicit list of the formal sentences that represent axioms in Th.

$\forall$x ( AXIOM(x) $\rightarrow$ PROOF(x) )

$\forall$x ( PROOF(x) $\wedge$ FOLLOWS_IN_ONE_STEP(x,y) $\rightarrow$ PROOF(c(x,y)))

$\forall$P ( ( $\forall$x AXIOM(x) $\rightarrow$ P(x) )
   $\wedge$ ( $\forall$x,y PROOF(x) $\wedge$ FOLLOW_IN_ONE_STEP (x,y) $\rightarrow$ P(c(x,y)) ) )
   $\rightarrow$ ( $\forall$x PROOF(x) $\rightarrow$ P(x) )

The next statement is not really a lemma. And although similar statements have appeared in the literature as lemmas I prefer to call it an observation.

> **Observation:** The relations described above define the relations they are named for.
> "Proof": Almost from the fact that $\mathbb{N}$ is formal. The formality of $\mathbb{N}$ means that all the operations required for theorem proving can be defined solely in terms of string manipulations – joining and splitting strings, that is concatenation and its inverse.
> Yet ultimately the comparison is to something extraformal. We have to argue that the relations based on concatenation are a reasonable representation of an extraformal notion, that they capture what meta-properties are conventionally agreed to be.
> Traditionally the debt to preformal ideas was expressed by saying

META defines a preformal string property Prop, if there exists a predicate P in META s.t.

For all a$\in\Sigma_{Th}$*, META |- P(a) iff Prop is true of a.

But this way of including unexamined properties is unsatisfactory. Concatenation-based string theories represent much the clearest expression, the least slippery grip we have on the concept of a string property. It would therefore in many ways be better to turn the definition talk around and say that the existence of a predicate in a string theory is what makes a well-defined string predicate. The predicates of string theories would then "define" by definition.

**Corollary:** The string relations defined above can be made to "correctly express".
Proof: Relations in a decidable theory that define express correctly. And all the string relations defined so far form part of a decidable theory, since

> By Gödel's first assumption, CONCAT exists and is consistent;
>
> By Gödel's second assumption, CONCAT extended by the axioms necessary to define the relations is consistent;
>
> By Gödel's third assumption, the consistent theory containing CONCAT and the relations can be further extended into a complete theory.

Label this theory META$_{Th}$.

> **Corollary:** META$_{Th}$ contains a predicate Prov() that expresses provability for the underlying theory.
> Proof: If PROOF(x) expresses correctly, then Prov(x) $\equiv \exists y$ PROOF(y) $\wedge$ ENDS(x,y) expresses provability for the underlying theory (with ENDS(x,y) defined as $\exists z$ c(z,x) = y).

Agreeing that Prov() correctly expresses is to agree that there is "nothing more to provability" than what can be laid down in formalisations in a concatenation-based theory.

**Proof of the Special Embedding Lemma:** Number theory is isomorphic to a subtheory of PROOF THEORY.

> **Definition:** The restriction of a theory to one of its predicates U(x) is the subtheory in which all predicates are of the form U(x) $\rightarrow$ P(x).
>
> U(x) is the universe of the subtheory.
>
> **Lemma:** PROOF THEORY restricted to NUMERAL(x) is isomorphic to number theory.
> Proof: The idea is to represent successor by concatenation with '. In this sprit we map

| the constant 0 | to | 0, |
| the successor relation s(x) = y | to | NUMERAL(x) $\wedge$ NUMERAL(y) $\rightarrow$ c(x,') = y, |
| the equality relation x=y | to | NUMERAL(x) $\wedge$ NUMERAL(y) $\rightarrow$ x=y. |

The axioms for concatenation in PROOF THEORY were chosen to make the Peano axioms provable under this transformation:

| 0 is a number <br><br><br> PROOF THEORY $\vdash$ NUMERAL(0) | This is the base clause defining NUMERAL(). |
|---|---|
| 0 has no predecessor <br><br> PROOF THEORY $\vdash \neg\exists x$ c(x;') = 0 | Immediate from the fact that 0 is a LETTER, and the primitiveness of letters, i.e. the axiom <br><br> LETTER(x) $\rightarrow \neg\exists u,v$ c(u,v) = x |
| Every number has a successor <br><br><br> PROOF THEORY $\vdash \forall x \exists y$ c(x,') = y | Immediate from the fact that c is a (total) function. |
| Two numbers are the same if they have the same successor.<br><br>PROOF THEORY \|- c(x,') = c(y,') → x=y | Follows from the second concatenation axiom in PROOF THEORY. |
| Induction<br><br>For all predicates P<br><br>PROOF THEORY \|- ∀P ( P(0)<br><br>    ∧ ( ∀x Numeral(x) → ( P(x) → P(c(x,')) ) )<br><br>       → (∀x Numeral(x) → P(x) ) | Induction in this or a logically equivalent form is an axiom of PROOF THEORY. |

Once successor is established, recursive definition gives addition and multiplication.

# Restructuring Cantor's Uncountability Argument

**N.B.** As models under the new definition no longer attempt to extensionalise, i.e. "cash in strings for objects", they do not have to provide a name string for every individual of the theory they are modelling. In fact, although they could voluntarily supply more names, they are constrained to provide names only for the (finitely many) constants of the theory. The idea of an uncountable model that contains more individuals than there could ever be strings does therefore not pose a threat: The new definition is able to deal with theories independently of whether their model in set theory would use a countable or an uncountable universe. It follows that the changes to postulates enabled in this section are not mandatory for the programme. The intention is only to illustrate the general way in which diagonal arguments may be restructured.

Cantor's diagonal argument forces a choice between two postulates. But since the traditional postulate has usually remained tacit some work needs to be done before the alternatives can be laid out.

> **Definition:** A function $p(n)$ from the numerals to the predicate strings of number theory that can be defined recursively in terms of concatenation is called a predicate sequence.

As a concession to standard usage write $P_n(x)$ for a sequence of predicates in a vector of variables x.

> **Definition:** A predicate sequence is said to be Cauchy convergent if
> for all M, for all n,m ≥ M, for all x ≤ M, $\mathbb{N}$ |- $P_n(x) \leftrightarrow P_m(x)$.

Predicate sequences that converge according to this definition will also converge (uniformly) as sequences of real numbers if the predicates are associated with a binary expansion valuing $P(n)$, the nth place, at $1/2^n$.

> **Definition:** A predicate L is said to be a limit for a predicate sequence if
> for all M, for all x ≤ M, $\mathbb{N}$ |- $L(x) \leftrightarrow P_M(x)$.

The minimal theory in which Cantor's argument can be made is now given by

> **Definition:** A predicate sequence space is a theory that
> a)  contains a model for number theory
> b)  allows concatenation-based predicate sequences to be defined

Working in a predicate sequence space, we can, for any small finite number n of predicate strings, construct by concatenation the predicate string:

> $D_n(x) \equiv$ "x=0 $\to \neg P_0(x) \wedge$ x=1 $\to \neg P_1(x) \wedge ... \wedge$ x=n $\to \neg P_n(x)$" (also written as: $\bigwedge_{i<n}$ x=i $\to \neg P_i(x)$ )

(Numberings and '…' are not symbols of standard number theory, so $D_n$ is not immediately a well-formed predicate string. It is what I would like to call a predicate expression – an expression that ostensibly denotes a predicate. When n is finite and small, however, it seems safe to assume that any such expression could be related back to a well-formed predicate string and therefore, since all predicate strings from number theory represent valid predicates, to a predicate.)

> **Definition:** Predicates in the sequence space are said to be standard (non-standard) if they are equivalent (not equivalent) to a predicate string already in available number theory.

Knowing the construction to be successful for small finite numbers of strings, one assumes that the expression denotes valid predicates also for large finite numbers. (To formally justify this conclusion, we would have to show the relative consistency of introducing into number theory two new elements of notation – numbering subscripts for predicates and finite quantification mixing subscript with term positions. By the informal argument already made, all such expressions can be eliminated in favour of standard predicate strings. Number theory extended by the new notations is not in fact any richer.)

Finally we come to the predicate expression "lim n→∞ $D_n$". The expression is clearly diagonal to all predicate strings covered by the numbering. (More formally: Limit signs cannot be eliminated in the same way as the first two additions to notation if we assume that the numbering includes all standard predicates.)

So conceding that the limit expression denotes a predicate is tantamount to the accepting a form of unformalisability: That there are predicates over the natural numbers not represented by predicate strings of number theory. But why are we obliged to accept the diagonal predicate for real? Could the implication of unformalisability not also count as a reason for rejecting the existence of a limit for the diagonal sequence?

At the decision for which types of sequences limits exist the metaphor of construction ends, infinitary postulates begin. Many different axioms of existence are conceivable to settle the matter, but two stand out.

The traditional favourite:

**Completeness Postulate –** The space contains a limit for every Cauchy convergent predicate sequence.

The alternative:

**The Church-Turing Axiom** – All predicates are standard.

which implies

**Conservative Convergence** – Only predicate sequences that converge to standard predicates have limits.

For most underlying theories, including number theory, the limits postulated by completeness are the maximum of what is consistently possible, conservative convergence describes the minimum.

The two postulates are incompatible according to

**Definition:** An extension of number theory is said to be productive (conservative) if it includes (does not include) non-standard predicates.

**Cantor's Major Theorem:** Predicate sequence spaces with the Completeness Postulate are productive.

The theorem follows directly from

**Cantor's Lemma:** There are predicate sequences that although Cauchy convergent do not converge to a standard predicate.
Proof: $(D_n)$ is one example of such a sequence. For details see Section Appendix.

For spaces that recognise only standard predicates a one-to-one correspondence between the set of natural numbers and its powers set is easy to establish. Any enumeration of the ordinary predicate strings of number theory will serve.

The predicate space under conservative convergence is thus consistent with the assumption that $|P(\omega)| = |\omega|$. Under the completeness postulate, which acts in the predicate space similar to the way that the Zermelo-Fraenkel axiom of comprehension acts in set theory, the cardinality of the power set must be strictly larger.

Insofar as both extensions lead to consistent theories, the Church-Turing Axiom and its denials have been shown to be independent of the axioms of number theory – independent of the second order axioms no less than first order.

Uncountability is thus something one postulates, not proves: The diagonal construction only unpacks the consequences of a postulate it does not attempt to justify. How many levels of infinity there is comparable to the familiar question in geometry of how many parallels run through a point. There is more than one plausible postulate, and therefore always a choice.

Even when the completeness postulate is rejected, the limit for the diagonal sequence denied to exist, Cantor's argument remains to prove a result that can appear to grow less trivial the longer it is contemplated, namely

**Cantor's Minor Theorem:** There are infinitely many predicates in number theory.
Proof: Suppose there were only a finite number. List representative strings for each equivalence class. String them together to produce $D_n$. The result is a valid predicate string that is provably non-equivalent to any of the strings in the initial list. Contradiction.

## Section Appendix

To illustrate the strength of the Completeness Postulate we will prove a more general result – the existence of a universal predicate from which the existence of Cantor's diagonal predicate follows.

Define a sequence t() from the numerals of $\mathbb{N}$ into the binary predicates of $\mathbb{N}$. The strictly literal definition stands on the right, a syntactically more relaxed version on the left.

$T_0 \equiv y = 0 \to P_0$                    $[ = c(\ \mathbf{y = 0} \to ;\ p(\mathbf{0})\ )\ ]$

$T_{n+1} \equiv T_n \wedge y = n+1 \to P_{n+1}$        $[ = c(\ t(n);\ \wedge\ \mathbf{y} =;\ c(n;\text{'});\ \to;\ p(c(n;\text{'}))\ )\ ]$

$T_n$ describes a growing table of predicates pushing out row by row from the origin.

It is evident from the definition of $T_n$ that the sequence is Cauchy convergent. For the remainder of this section, 'convergence' will be used as a shorthand for 'convergence to a standard predicate'.

**Cantor's Lemma:** The limit of $T_n$, if it exists, is a universal predicate. This is to say that if T is the limit of $T_n$, and $g=p^{-1}$, then for all predicate strings P, numerals n, m,

$$\mathbb{N} \vdash T(g(P),m) \text{ iff } \mathbb{N} \vdash P_n(m).$$

Proof: Assume that the limit exists, and write T for it.

Pick any n,m, and choose M larger than both.

Then from the definition of T, $\mathbb{N} \vdash \forall i,j < M\ T(i,j) \leftrightarrow T_M(i,j)$.

Hence $\mathbb{N} \vdash T(n,m) \leftrightarrow T_M(n,m)$.

We are done if we can show that $\mathbb{N} \vdash T_M(n,m) \leftrightarrow P_n(m)$. But this is a propositional tautology of the

form:

$$( (\perp \to P_1(n)\ ) \wedge ... \wedge (\perp \to P_i(n)\ ) \wedge ... \wedge (\ T \to Q\ )\ ) \leftrightarrow Q,$$

which $\mathbb{N}$ will prove.

**Corollary:** The diagonal predicate sequence defined by $d(n) = \neg T_n(x,x)$ does not converge

Proof: Let P be any predicate string in $\mathbb{N}$, let $n = p^{-1}(P)$.

And suppose that d() were to converge to P, i.e. $\forall M\ \mathbb{N} \vdash \forall i \leq M\ P(i) \leftrightarrow D_M(i)$.

Then in particular, $\mathbb{N} \vdash \forall i \leq n\ P(i) \leftrightarrow D_n(i)$. And again, $\mathbb{N} \vdash P(n) \leftrightarrow D_n(n)$.

Now p(n) = $P_n$. Therefore, relabelling only, $\mathbb{N}$ |- $P_n(n) \leftrightarrow D_n(n)$.

On the other hand, by the definition of d(), $D_n(n) = \neg T_n(n,n) = \neg T_n(g(P),n)$.

Because of string identity, we have $\mathbb{N}$ |- $\neg T_n(g(P),n) \leftrightarrow D_n(n)$ from rewriting the (provable) tautology

$D_n(n) \leftrightarrow D_n(n)$.

Now from the theorem, $\mathbb{N}$ |- $\neg T_n(g(P),n) \leftrightarrow \neg P_n(n)$.

Hence $\mathbb{N}$ |- $\neg P_n(n) \leftrightarrow D_n(n)$.

Contradiction, as was first shown by Cantor over a hundred years ago.

**Proposition:** $T_n$ converges iff $\neg T_n$ converges.

Proof: $\mathbb{N}$ |- ( $\forall i \leq M$   $T(\mathbf{x}) \leftrightarrow P_i(\mathbf{x})$ ) $\leftrightarrow$ ( $\forall i \leq M$   $\neg T(\mathbf{x}) \leftrightarrow \neg P_i(\mathbf{x})$ ), a propositional tautology.

**Proposition:** If the diagonal of a binary predicate diverges, the binary diverges.

Proof: The existence of an M s.t M does not prove $\forall i \leq M$   $P(i,i) \leftrightarrow P_i(i)$ readily implies that $\mathbb{N}$ does not prove the more general $\forall i,j \leq M$   $P(i,j) \leftrightarrow P_i(j)$.

**Corollary**: The sequence $T_n$ does not converge.

# Restructuring Set Theory

**N.B.** As the new definition of a model no longer relies on set theory, the postulate changes in this section are not strictly necessary for success. However, the programme is highly sympathetic to the generalisation of comprehension attempted here. The main point to take away: There is no one catholic theory of sets that could serve as a foundation for mathematics, but a multitude of widely different versions.

The restructuring of set theory begins with a choice of two postulates.

> **The traditional favourite:** The Zermelo-Fraenkel axiom of comprehension
> $$\forall P \ \forall a \ \exists p \ \forall x \ \ x \in p \leftrightarrow P(x) \wedge x \in a$$

> **The alternative:** V, the set of all sets, exists.

A diagonal argument to show that the two postulates are incompatible:

> The two assumptions combined are equivalent to the naïve axiom of comprehension,
> $\forall P \ \exists p \ \forall x \ x \in p \leftrightarrow P(x)$, which Russell's predicate $x \notin x$ shows to be inconsistent.

Three reasons for believing that the alternative is at least as good:

1) The alternative is the more conservative and has great intuitive appeal. The Zermelo-Fraenkel axiom, on the other hand, is an example of that most undesirable sort of assumption: strong yet intransparent.

2) Only when V exists can set theory extend the Boolean algebra. Under the Zermelo-Fraenkel axiom no set has a complement that is also a set. So although the Boolean operators of union and intersection are available, the third operator is missing. This produces an unattractive loss of symmetry in the Zermelo-Fraenkel system, as structure growing up from $\varnothing$ is not matched by structure growing down from V.

3) The Zermelo-Fraenkel axiom makes distinctions between infinite diagonalisations – banning some (e.g. Russell's paradox) while allowing others (e.g. Cantor's uncountability argument) – that are hard to defend on principled grounds.
   Since every set is a subset of V the power set of V equals V. The identical function i from V to V is therefore a bijection between V and P(V). Given that $\{ x \mid x \notin i(x) \} = \{ x \mid x \notin x \}$, Russell's diagonal argument applied to V and i runs exactly parallel to Cantor's diagonal argument applied to a bijection between $\omega$ and $P(\omega)$.

A thesis ventured in conclusion:

> Compared to the alternative the Zermelo-Fraenkel axiom of comprehension could seem too narrow because it excludes V and the identical function from V to V, too wide because it admits the diagonalisation of infinite sets. In others words, there is more than one viable variety of set theory. Comprehension needs to be generalised to encompass the full range of varieties.

The task which any replacement for the familiar axiom of comprehension will have to take over is sorting the predicates of set theory into those that can be extensionalised, and those that cannot. While it is easy to give upper and lower bounds tracing the exact boundaries is going to require a shift in perspective.

Write EXTEN(P) for $\exists p \ \forall x \ x \in p \leftrightarrow P(x)$ with varying P. The sentence says that the predicate P has an extension, the set p.

The upper bound for comprehension is the naïve axiom, $\forall P \ \text{EXTEN}(P)$.

Comprehension confers on predicates to ability to serve as names for sets. In this capacity, they can, like any name, turn out to be empty.

> **Definition:** A predicate that cannot be extensionalised ($\neg$EXTEN(P) is provable) is called empty.

**Russell's Major Theorem:** There are empty predicates.
Proof: Use Russell's predicate x∉x.

The case of Russell's predicate is easy. It can never be consistently extensionalised. But there are other predicates that might plausibly be extensionalisable in some versions of set theory, and not in others. For example, the constantly true predicates standing for V, the set of all sets. If despite this evidence we insist on a single unified theory of sets, the status of these extensions is rendered ambiguous. Quantifying over versions produces a clear result only for predicates that are extensionalisable in all versions ('sets'), and impossible predicates like Russell's that are extensionalisable in no version. It leaves the extension of predicates that are extensionalisable only in some versions ('classes') in ontological limbo.

Traditionally, classes were conceded to exist, directly or indirectly; but not without awkward restrictions. Provided the theory lost no power any designer of axioms would clearly prefer to introduce only one type of extension – sets with full operational privileges –, and to be able, for every predicate, either to prove or refute the existence of its extension. By allowing for more than one version of set theory this may indeed be possible.

For every particular version of set theory the extensionalisability of a predicate is fully determined, even though it may of course vary across versions. Empty predicates no longer have bizarrely behaved, unsetlike extensions (e.g. the shadow classes of the Zermelo-Fraenkel system, or the explicit classes of Neumann-Bernays-Gödel theory), they have no extension at all. Rather like limitless sequences, empty predicates are second order strings that do not convert into a first order object. They exist as names, and names only.

A lower bound on comprehension is the axiom restricted to finite domains,

$$\forall P \ \forall finite \ a \ \exists p \ \forall x \ \ x \in p \leftrightarrow P(x) \land x \in a$$

(with *finite* defined in the usual sense that no bijection exists between the set and a proper subset)

Only over finite domains is the intuition that suggested naïve comprehension fully trustworthy. Each instance of comprehension over infinite sets amounts to an infinitary postulate that can only be granted after careful justification.

Finite comprehension is linked to an extension of the Boolean algebra by the successor function, s(x) = {x}.

**Fact:** Equations in the six functions ∅, V, ∪( ), ∩( ), ¬( ), s( ) are isomorphic to predicate sentences in ∈. (It makes little difference whether one treats ∅ and V as constants or constant functions).
Proof (Abbreviated): The critical step in the conversion is x∈y ⟷ s(x) ∩ y = s(x).

This theory, which can be expressed with only functions and equations, or predicates and sentences, or a mixture of both, might be called the core of set theory.

**Proposition:** Finite comprehension is implied by the core.
Proof: The existence of all finite subsets of a given sets is provable from s( ) and ∪( ).

In the context of finite comprehension Russell's diagonalisation is revealed as a paradox of infinity.

**Russell's Minor Theorem:** There are infinitely many sets.
Proof: Run Russell's predicate against the finite comprehension axiom.

Over the core, the axiom of infinity thus becomes an axiom of arithmetic – no longer needed to introduce infinite sets but ω, the structured set of the natural numbers.

**Fact:** The core of set theory cannot be satisfied over a finite domain. The largest subtheory that can is known as the theory of a field of sets.
Proof (Abbreviated): Dual theorems involving s( ) extend the Boolean algebra to the field of sets, dual and non-dual theorems together extend it to the core of set theory.

Between these upper and lower bounds lie candidate axioms, both skewed ones like the Zermelo-Fraenkel axiom and symmetric expansions of the core. In order to make set theory as powerful as we can, we would like to see the candidates introduce as many extensions as possible.

**Definition:** A shot at set theory is a maximally consistent set of sentences containing, apart from the axiom of extensionality, only sentences of the form EXTEN(P).

It is tempting to think that set theory could be summarised by the following

**Meta Axiom:** The axioms of any set theory form a shot.

The idea is to give, instead of explicit lists of axioms for each version of set theory, a general criterion the axioms for any version must satisfy. Characterising set theory in this way as the solution to an equation brings freedom almost as great as naïve comprehension. It encourages a hypothetical style of reasoning that draws on whichever existence assumptions seem necessary for a derivation. Only when a contradiction arises, showing that two sets of assumptions of the form EXTEN(P) are mutually inconsistent, is the need for a decision established. At that point set theory branches into two systems. By a not very scientific process usually one of these systems (the one that "feels" more natural) will be designated the standard system, the other as non-standard. So instead of blindly guessing a system – and none of the axiom systems for set theory proposed so far have been satisfactory – one could begin to explore the shape of possible solutions, slowly to narrow down their range.

**Definition:** A sentence is true for set theory if it is provable in all shots.

Russell's Major Theorem is one of the relatively few sentences true for set theory. The large remainder will be provable only in more restricted shots – each a version of set theory. A first divide is between those versions that are Boolean, and those that are not.

**Definition:** A shot is Boolean if it is maximally consistent subject to these constraints:
$$\text{EXTEN}(P) \Rightarrow \text{EXTEN}(\neg P)$$
$$\text{EXTEN}(P) \wedge \text{EXTEN}(Q) \Rightarrow \text{EXTEN}(P \wedge Q)$$

Note that because of reverse comprehension, $\forall p\ \exists P\ \forall x\ P(x) \leftrightarrow x \in p$, Boolean shots are consistent with the existence of the three Boolean functions.

**Proposition:** In a Boolean shot, $\neg \text{EXTEN}(x \in x)$.

A Boolean shot excludes non-complemented predicates, which exhibit many of the strange properties associated with the non-recursive sets of recursion theory.

For n-ary predicates, let the extension p consist of ordered n-tuples.

**Proposition:** In a Boolean shot, $\neg \text{EXTEN}(x \in y) \vee \neg \text{EXTEN}(x = y)$.
Proof: $x \in x$ is the conjunction of $x \in y$ and $x = y$.

That the main relations of set theory fail to be extensionalisable perhaps explains how classes were invented to preserve an extensionalisabitity of sorts for set theory.

**Definition:** A grape is a sentence from a shot.

As we have already seen in the Boolean case, one way of systematising shots is to postulate the existence of functions, which effectively condenses infinitely many grapes into a single assumption.

The next obvious functions after the Boolean trio are successor and the power set function. More ambitiously, one could postulate for some shots the existence of all recursive functions over given sets.

**Definition:** A shot is recursive if it is consistent with the existence of all recursive functions.

**Proposition:** There are non-recursive functions, i.e. function strings that need not be extensionalisable even in recursive shots.
Proof: Consider the string describing the indicator function of Cantor's diagonal set.

The extensionalisability assumption for Russell's predicate is inconsistent by itself. Excluding Cantor's predicate from shots requires a further axiom.

**The Axiom of Single Infinity:** "There exists a bijection between any two infinite sets"

Single infinity is consistent with the existence of V. It replaces assumptions such as the axiom of choice or the continuum hypothesis that are independent of the elementary Zermelo-Fraenkel system.

The axiom states – truthfully – that being infinite is still compatible with almost any structure. Knowing a set to be infinite is to know very little about it.

It may violate intuitions schooled on Cantorian precepts, but the axiom is no more contrary to common sense than the generally accepted fact that the unit interval contains the same number of points as the real line, the real plane, indeed as any finite-dimensional real space.

**Definition:** A standard shot is a Boolean shot that is also consistent with the axiom of single infinity.

Because of diagonal arguments, every shot has to make decisions between bijections and diagonalisations. Standard shots consistently decide in favour of the former.

The axiom of single infinity contained in any standard shot makes it possible to refute existence assumptions leading to uncountability. There are many, and they can be intricately hidden. Only by dropping the principle of separation on which the Zermelo-Fraenkel axiom is based – EXTEN(P) $\Rightarrow$ EXTEN(P$\wedge$Q), for any Q – can one hope to become aware of these subtleties.

**Proposition:** EXTEN(P) and Q$\rightarrow$P do *not* imply EXTEN(Q) (even though Q(x) $\equiv$ Q(x) $\wedge$ x$\in$p).
Proof: A counterexample is given by Russell's predicate and x=x.

The fact that Russell's predicate, R, implies any constantly true predicate only shows hypothetically that if R did have an extension it would be a subset of V. The isomorphism between the $\rightarrow$ relation over predicates and the $\subseteq$ relation over sets holds only for extensionalisable predicates.

The point of set theory is to discover, not declare, which predicates have extensions.

## Section Appendix

Mathematical imagination depends on pictures at least as much as on formal derivation. For new concepts to take root mental pictures also have to change.

The picture below wants to argue that Russell's diagonalisation is essentially the same as Cantor's, and that the success of both constructions is evident only over finite domains.

| $\in$ | $S_1$ | $S_2$ | $S_3$ | ... | $S_{N-1}$ | $S_N$ | $S_{N+1}$ | ... |
|---|---|---|---|---|---|---|---|---|
| $S_1$ | | ▓ | | | ▓ | | | |
| $S_2$ | ▓ | ▓ | ▓ | ... | | ▓ | | |
| $S_3$ | | ▓ | | | | | ▓ | ... |
| ... | | ... | | ... | | ... | | |
| $S_{N-1}$ | | | ▓ | | ▓ | ▓ | ▓ | ... |
| $S_N$ | | ▓ | ▓ | | | | | |
| ... | | ... | | ... | | ... | | ... |

The meaning of any picture of diagonalisation turns on what is inevitably only hinted at with three dots or a similar device: Does the table, does its diagonal actually reach infinitely far ? Or does it come up against a flexible barrier, a barrier that allows the table to reach arbitrarily finitely far but no further ?

Diagonal arguments effectively postulate the existence of an infinite table, so that they can then prove the existence of the table's diagonal. Once the postulate is made, the existence of the diagonal is self-evident. But the postulate itself is not self-evident. The existence of the whole table is at most heuristically suggested as one of several possible postulates; in the case of set theory it is even inconsistent.

Over any finite domain, the membership relation could be specified extensionally. Given that it is not unrealistic to list all sets – simply by listing the predicate strings of set theory, as theories are restricted to the objects their predicates can distinguish – the limit of these approximations would represent everything there is to know about membership, and hence, set theory.

$$
\begin{array}{lll}
S_1 \in S_1 & S_1 \notin S_2 & S_1 \in S_3 \quad \cdots \\
S_2 \notin S_1 & S_2 \in S_2 & S_2 \in S_3 \\
S_3 \in S_1 & S_3 \notin S_2 & S_3 \notin S_3 \\
\quad \cdots & & \quad \cdots
\end{array}
$$

And for this reason the limit, the actually infinite table, is never reached.

# Restructuring Computability

**N.B.** The rejection of the traditional postulate isolated in this section is mandatory for the success of the programme. This is understandable since the changes made here retrace those explained in much greater detail in the Gödel section.

## I.    Recursion Theory

The predicate calculus is the gold standard for predicate expressions. Recursion theory often uses more idiosyncratic notation, but fortunately almost all can be related back to the predicate strings of either number theory or PROOF THEORY.

For a formalism that allows all recursive means of definition write

> S for the class of all predicate expressions in terms of successor,
> C for the class of all predicate expressions in terms of concatenation (over an alphabet $\Sigma$).

We know that S is isomorphically included in C – the ability of recursively enumerating sets of strings includes the ability to enumerate sets of numerals.

C, however, is much more extensive than S, and its relationship to S more complex.

Defining predicate expressions can do two things. It can describe which expressions are well-formed. And it can deliver a stamp of approval – declaring which expressions are genuine predicates and "really" represent a subset of the universe. The first is innocuous, the second effectively hides an existence assumption.

> **Definition:** A predicate expression P from a given formalism is said to be solid (empty) if it there exists (does not exist) a model in set theory in which m(P) is extensionalisable.

Put more simply, a predicate is solid if can be consistent to assume that it has an extension. To motivate the definition, recall how Russell's paradox essentially rehearses the fact that excluded middle breaks down when empty names are introduced into a logic: The apparent tautology Bald(x) $\vee$ ¬Bald(x) is false when applied to 'The present king of France' as much as x$\in$x $\vee$ x$\notin$x is false when applied to 'The extension of the predicate x$\notin$x' (={x | x$\notin$x }).

All expressions in terms of successor are solid. The limits of non-convergent predicate sequences are examples of empty predicate expressions. Like the reigning roi, who is neither bald nor hirsute, the extension of empty predicate expressions is neither identical nor non-identical to solid predicates.

Recursion theory now has to make a choice. At most one of the following two assumptions can be declared true by definition.

The traditional assumption has been that all predicate expressions give rise to a set.

> **Naïve Extensionalisability:** All predicate expressions in C are solid.

The naïve extensionalisability assumption can be regarded as the string-theoretic equivalent of the naïve comprehension axiom for set-theoretic predicates.

The alternative is to concede fewer extensions. While there are millions of intermediate possibilities, the most natural among the alternatives is also the most restrictive. It consists of tying the concept of a set of strings to recursiveness. This move in effect elevates the Church-Turing Thesis into an axiom[§].

---

[§] One way of arriving at the Church-Turing Axiom is through a syllogism:

<div>

Church-Turing Thesis: "All formalisable functions are recursive"
+ Formalism: "All functions are formalisable"

Church-Turing Axiom: "All functions are recursive"

</div>

**Church-Turing Axiom:** Sets of strings are only well-defined if they are recursive (i.e. both the set and its complement are represented by expressions in C).

Naïve extensionalisability says that all predicate expressions are associated with a set; the Church-Turing Axiom says, conversely, that all sets of strings are associated with a predicate expression.

The two postulates are incompatible according to

**Church-Turing Theorem:** There are predicate expressions in C that, if assumed solid, would give rise to non-recursive sets of strings.
Proof: As provability and substitution are definable in terms of concatenation, C can contain a universal predicate for any formalism. The diagonal predicate for C's own formalism, if it had an extension, must then be non-recursive.

A definition of the form "$d(n) = f_n(n) + 1$" presupposes the existence of a universal predicate, $u(n,m) = f_n(m)$. That a universal predicate could be "extensionally true" – behave as intimated by the crossover of n from subscript in $f_n(m)$ to argument in $u(n,m)$ – is an assumption one can make, but it is hardly self-evident, and has the even less evident consequence of undecidability.

The consequences of postulating undecidability are well-known, and need not be rehearsed again. For the remainder of this section we focus on theories that agree on the need to deny naïve extensionalisability.

Diagonalisation is inherently non-recursive, whether introduced into recursion theory by unconventional notation or into set theory by the axiom of separation. There are hence three types of string theories: Standard theories in which the Church-Turing Axiom holds; non-standard theories that admit more than the standard extensions and/or non-recursive construction methods but remain consistent with decidability; and non-theories that pretend to prove undecidability.

In the standard theories, once non-recursive sets are defined out of existence by the Church-Turing Axiom, proving the existence of a non-recursive set becomes tantamount to proving a contradiction. Under the alternative axiom, the Church-Turing Theorem is thus employed to refute the idea that all predicates in C are solid. So rather than assuming the existence of certain extensions and using it to refute decidability, we can also assume decidability as basic and use it to refute the existence of certain extensions.

The natural assumption still remains that a predicate expression should be treated as solid unless the existence of the set it would give rise to can be shown to imply undecidability (all theories) or non-recursiveness (only standard theories). We can safely assume that finite fragments of extensions exist; what we cannot assume is that the fragments will always assemble into a whole.

**Definition:** A finite fragment of a predicate expression in C is the extension of $P(x) \land x \in S$, for any finite set of strings S.

Let P be an empty predicate expression. Write $|x|$ for the string length of x. For any n, $P_n(x) \equiv P(x) \land |x| < n$ is solid. The set of all strings of length less than n such that P is recognised to exist. Saying that P itself is not solid is saying that $\lim n \rightarrow \infty P_n$ does not exist.

Empty expressions do not give rise to non-recursive sets. They give rise to no set at all. They also have no equivalent among the successor-based expressions under a translation.

A sharp distinction must be drawn between predicate expressions in S, all of which have extensions, and attempts to link these predicates to expressions in C through a code. The link fails because predicate expressions at the other end of the link fail to be extensionalisable. For the string sets ostensibly named by expressions in C, 'recursively enumerable' should be clarified to mean 'certain to exist finitely far, limit only if proven safe'.

**Proposition:** Extensionalisability is not preserved by composition.
Proof: The extension of the monadic predicate Prov() is recursive and exists for any decidable theory. Substitution is primitive recursive, so its extension ought to be granted even more readily. The dyadic predicate Prov(sub(,)), however, cannot be extensionalised.

"The universal predicate" is an empty name rather like "The extension of the membership relation $x \in y$".

## II.  Turing machines

Turing machines are a type of successor-based predicate expression.

> **Definition:** The typographically represented set of instructions of a Turing machine is called a Turing programme.

The function a Turing machine computes should be identified with its (second order) Turing programme, not the (first order) movements on tape. The Turing programme represents a completed infinity, which is a much better way of thinking about functions than to imagine an open-ended process of computation. A process is by its nature ambivalent between 'for any finite number of steps/inputs' and 'for an infinite number'.

Successor may be granted, but can Turing machines compute concatenation ?

The head of a Turing machine scans only one square of its tape. The machine has no intentions, and should not be credited with being able to see beyond the boundaries of the square it presently scans.

Until we see the machine respond to an input of

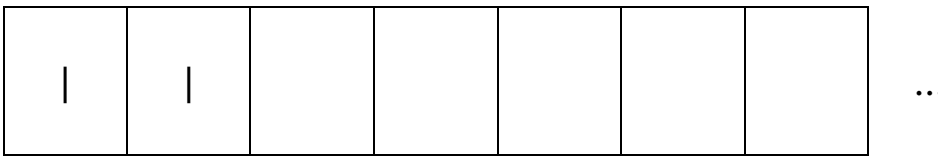

with an output of

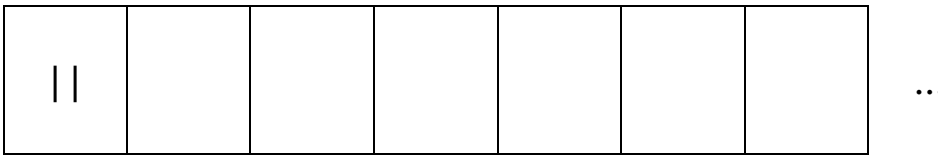

we would have to conclude that the concatenation of inputs across the boundaries of fields was an involuntary addition made by the observer.

The interpreter, whether human or another machine, who performs the concatenation to find that an input of strokes stretching over several squares represents evidence of the machine having accepted an input of their concatenation must be able to recognise many more symbols than the machine being interpreted. The concatenating machine is thus a type of meta machine.

As inputs into the underlying machine stretch over more and more squares of the tape, the interpreting machine has to recognise and put out more and more symbols on its own squares to stay ahead. At the limit, the interpreting machine gets squeezed. A limiting machine would have to be able to recognise infinitely many different inputs of symbols on a single square of its tape – which violates the definition of a Turing Machine, and any sensible definition of computability as well.

The assumption that Turing machines can compute general sets of strings is equivalent to the existence of CONCATENATOR, a machine that, before a Turing machine TM starts, deconcatenates (spreads) the input in the leftmost field over the appropriate number of squares on the tape and after TM halts concatenates the contents of all non-blank fields back into the leftmost field.

> **Turing's Theorem:** CONCATENATOR does not exist.
> Proof: CONCATENATOR would be able to do what META cannot.

The only possible treatment of concatenation is extensional, which conflicts with the inability of agents to handle infinitely many objects.

# Inventing a Language for Paradox

An alternative way of summing up the dilemma between the incompleteness or inconsistency of META is to say that concatenation and the predicate calculus do not mix.

So if the predicate calculus, or other standard formulations like the Turing machine, fail at the limit to contain concatenation then there is a case to be made for looking at an altogether different kind of language.

As PROOF THEORY$_L$ is in effect the instruction manual for handling languages L these revisions take the form of weakening or strengthening the axioms of PROOF THEORY$_L$.

> **Definition:** The predicates LOGICAL_CONSTANT, VARIABLE, OPEN_TERM, CLOSED_TERM, PREDICATE, SENTENCE in PROOF THEORY are said to be the roles of strings.

> **Definition:** A string X is said to be cast into a role ROLE if PROOF THEORY$_L$ proves ROLE(X).

The predicate calculus is distinguished by a set of strict exclusion clauses.

> **Definition:** The Union Regulation Axiom in PROOF THEORY$_L$ states that for every pairing of roles ROLE$_i$, ROLE$_j$, i≠j
> $$ROLE_i \rightarrow \neg ROLE_j(x) \wedge ROLE_j(x) \rightarrow \neg ROLE_j$$

The regulation implies no obligation to hire. A string could be cast into neither ROLE$_i$ nor ROLE$_j$ and there are in fact hopeless strings ("noise") that are not cast for any of the six roles. What is excluded is that a string actor play more than one role.

We now define a type of theory that breaks radically with the conventions of the predicate calculus.

> **Definition:** A language L is said to be wild if PROOF THEORY$_L$ proves
>
> $$\forall x\ OPEN\_TERM(x) \leftrightarrow PREDICATE(x)$$
>
> $$\forall x\ CLOSED\_TERM(x) \leftrightarrow SENTENCE(x)$$

> **Definition:** A circular language is a wild language where PROOF THEORY$_L$ also proves
>
> $$\forall x\ OPEN(x) \leftrightarrow CLOSED(x)$$

In a circular language strings can ambiguously be both, open and closed, predicates and functions.

> **Theorem:** Circular languages contain a concatenation function in the sense that
> $$\exists c\ \forall xy\ c(x,y) = xy.$$

> Proof: The legally available functional xy already is, to all intents and purposes, concatenation. And the language obviously proves $\forall xy\ xy = xy$. Then merely relabelling xy on one side of the equality proves the existence of c.

Among ordinary languages the sequence CONCAT$_n$ has no limit; the circular language, insofar as it is able to derive more than finitely many concatenation statements, could be spoken of as a limit. Circular languages are the real home of (infinite) concatenation.

Legalistic, but formally unobjectionable:

> **Theorem:** Circular languages prove $\forall x,y\ x=y$.
> Proof: Since x≠xy, y≠xy, xy is a variable not bound by either $\forall x$ or $\forall y$. Any unbound variable may be replaced by any other free variable.
> So from $\forall x,y\ c(x,y) = xy$ we can derive $\forall xy\ c(x,y) = z$, assuming only axioms.
> Since z was free, this gives: $\forall x,y,z\ c(x,y) = z$.

**Corollary ("Out of Eden"):** Circular languages are only consistent with a universe of size one or less.
Proof: Take any two primitive terms a,b, a≠b. Then c(a,b) = a, c(a,b) = b. Hence a=b. Contradiction. Therefore ∀x,y x=y.

The language breaks down as soon as it is expelled from the paradise of undivided oneness. Collapsing the distinctions between the roles of strings has the effect of merging the whole theory into one undifferentiated blob.

But this is not the end of the road. The result can still be improved on.

**Definition:** In a wild language, the merged role of PREDICATE ∨ OPEN_TERM is called FUNCTOR.

**Definition:** A wild language L is said to be a Truth-and-Substitution Language if there exist functors TRUE and SUB such that ∀X TRUE(X) ↔ X and ∀P∀X P(X) ↔ SUB(P,X) are theorems of L.

**Lemma:** A circular language that is consistent with a universe of size greater or equal one has a truth predicate.
Proof: Fed any term S, the language proves ∀S S ↔ ( S ↔ S ). As there is no principled distinction between sentences and closed terms we can think of S ↔ S as the predicate "$S_{↔}$(S)". So we have $S_{↔}$(S) ↔ S. Relabel $S_{↔}$ as TRUE, and you have a truth predicate.

**Lemma:** A circular language that is consistent with a universe of size greater or equal two has a substitution predicate.
Proof: Fed two terms P and a, the language proves ∀P ∀a P(a) = P(a). From rewriting the right hand side we have P(a) = Sub(P,a). Since P(a) and Sub(P,a) are equal as terms they are certainly equivalent as sentences. Therefore ∀P,a P(a) ↔ Sub(P,a).

**Corollary:** A circular language that is consistent with a universe of size greater or equal one has a diagonal substitution predicate SSB s.t. ∀P P(SSB) ↔ SSB(SSB).
Proof: Same argument as above with P=a.

**Theorem ("The General Paradox"):** Truth-and-Substitution languages are inconsistent.
Proof: Define FSSB := ¬TRUE(SUB(X,X)).Then given that FSSB(SSSB)) = ¬TRUE(FSSB(SSSB))),

$$\text{FSSB(FSSB)} ↔ ¬\text{FSSB(FSSB)}.$$

**Corollary:** Circular languages can be satisfied only by the empty universe.
Proof: Suppose ∃X. By the lemmas above the language can then define FSSB.

This is a strong result. It would be wrong, however, to dismiss circular languages as trivial because of it.

Circularity has advantages. Observe the naturalness of the "axiom of concatenation" ∀x,y c(x,y) = xy and the elegance and intuitiveness of TRUE(X) ↔ X. And contrast this with the puritanical awkwardness of the predicate calculus (inverted commas, corners, ...) when it comes to dealing with the same concepts. The insistence on enforcing a strict separation of types, between meta and base, around which the predicate calculus is built is legitimate but has a price.

As far as one can determine these things most natural languages are in fact circular, or stand at least much closer to the insouciance of wild languages than the rigour of the predicate calculus. But if natural languages are broadly circular (and suckers that we are, their users tend to believe in the existence of objects) how do they avoid inconsistency ? Well, everyday life works with a background assumption that there are only finitely many objects, or at least that one can only deal with finitely many objects at any one point in time. This can still give rise to a potential sort of infinity by shifting attention from one set of objects to another, and another.

Now over finite domains type violations, or pulling objects down from one layer of stratification to another, turn out to be harmless. Let the general case be inconsistent as every recurring paradox shows, but restricted to a movable finite domain wild languages can be perfectly valid.

**Definition:** A wild language from which we drop the solidity assumption (that the mere availability of a function or predicate letter implies existence) and then make the truth of relations and the existence of values for functions conditional on ∃N ∃$A_i$ ∀X ( X=$A_1$ ∨ X=$A_2$ ∨ .... ∨ X=$A_N$ ) is said to be tamed.

As usual in mathematics it makes no difference whether N is set to be a small number that can be effectively handled, or larger than the number of atoms in the galaxy. The salient distinction is between finite and infinite, not different sizes of the finite. The constraint, although real, need not figure as such for any and all practical purposes. People can happily live their lives immersed in wild languages and never encounter a problem (unless, that is, they are Greek and called Epimenides).

In a tamed language the paradox simply proves that the diagonal object is out of scope, and looking out from a small fenced-in patch there is certainly much scope outside of it. Refocusing to the patch next door is analogous to exchanging the pre-declared list for another.

**Definition:** A theory Th is said to be pseudo-circular if there exist in Th a predicate T and a binary function sub and outside of Th a function g: STRINGs → CLOSED_TERMs such that

For all sentences S, Th |- T(g(S)) ↔ S and

For all predicates P, terms t, Th |- sub(g(P),t) = g(SUBSTITUTION(g(P),t))

**Corollary**: In a pseudo-circular language there exists a predicate U ("the universal table") such that

For all predicates P, terms t, Th |- U(g(P),t) ↔ P(t)

**Corollary**: Pseudo-circular languages are inconsistent.
Proof: A pseudo-circular language contains the image of a Truth-and-Substitution-Language. By the General Paradox, the image is inconsistent. Therefore the language is.

**Lemma:** Set theory under naïve comprehension is pseudo-circular.
Proof: It is said that Peano, when he invented the ∈ notation, derived it from εστι, the Greek for "is". With that we have a∈P ⟺ a "is" P ⟺ P(a) ⟺ TRUE(P(a)). In other words, with g the function that sends predicate strings P to term strings {x | P(x)}, the reverse of membership, ϶, is a universal table (and Russell's Paradox yet another gloss on what we have called the General Paradox).

**Lemma:** META, if it existed, would be pseudo-circular.
Proof: For any theory Th powerful enough to host its meta theory, META_{Th} contains a truth predicate and a substitution function. The function g can be set to be the modelling function from META to Th.

**Corollary:** Tarskian model theory is pseudo-circular.
Proof: There are two ways of seeing this. First note that according to the traditional definition of a model that goes back to Tarski the predicate True_in_a_model is expected to do what only a provability predicate in a meta theory could. So the circularity of model theory follows from the circularity of META.
Second observe that in Tarskian model theory it is assumed that

(∀ P in Th) (∃ p in M) (∀ a in M)  P(a) ↔ |a|∈ p.

Where | | represents the "cash"-function that converts general closed terms into canonical names, e.g. (0'+0) into 0'. Provided Th contains canonical names for all its objects, the cash function can be eliminated and the universe of Th and M merged:

(∀P∈ Th) (∃p) (∀a∈u)  P(a) ↔ a∈ p.

But this is in effect the naïve axiom of comprehension. Hence the circularity of model theory also follows from the circularity of set theory.

**Lemma:** If the limit of Cantor's sequence T_n existed, number theory would be pseudo-circular.
Proof: See appendix to Cantor section.

**Grelling's Lemma:** Almost natural English is pseudo-circular.
Proof: Grelling's paradoxical property heterological() is the predicate diagonal to a universal table of all English predicates of the form is_X, where X is an adjective. P is heterological iff ¬P(P) (= Fssb(P)).

**Epimenides' Theorem:** Natural languages are pseudo-circular.
Proof: The simple reason is that natural languages contain truth, substitution, and concatenation predicates. For the long version we refer to a famous Cretan.

To assimilate the liar paradox − This sentence is false − to the general paradox for

This Sentence has property P

write

P(This_Sentence).

The value of This_Sentence(), construed as an operator on strings, could be paraphrased as

"The string This_Sentence substituted into the enclosing predicate context P."

Granted that This_Sentence() is a function of predicate strings, how could it be a function of its own context predicate ? By diagonalisation − equating P and Q in Q(func(P)) effectively turns a context-sensitive string function func(P) into a function of its own predicate context Q≡P.

It is not hard to recognise This_Sentence as ssb(). Under diagonalisation, the self-substitution function ssb inserted into P does precisely what one would expect of This_Sentence:

ssb(Pssb) = Pssb(Pssb), with the effect that Pssb(Pssb) is self-referential in P.

So putting This_Sentence(P) ≡ ssb(Pssb) and P≡Is_False, the reconstructed liar paradox

Is_False(This_Sentence(Is_False))

becomes another way of writing the general paradox

Fssb(Fssb).

# Restructuring the Rest

The present situation may run parallel to the historical choice between Euclidean and non-Euclidean geometries, except that up to now we seem to have been living in the non-standard, bent world. The standard world is defined by the axiom of single infinity, non-standard worlds by any form of denying single infinity.

Results that were proved based on assumptions that hold true only in the bent world will of course remain provable within in their own protected sphere. Outside in the big wide straight world, however, their survival is not guaranteed. Some can be adapted, others will evaporate.

The quickest way to isolate non-standard assumptions is to re-examine anything that could imply uncountability. To treat, for the time being, a proof of uncountability as a proof of $\perp$. (In set theory with the axiom of single infinity any derivation of uncountability automatically amounts to a proof by contradiction.)

Not all results mentioning higher cardinals will simply be lost. Even under single infinity some results that were originally stated as being about cardinality may be recovered as results about orderings or other structure inherent in the set. Recovery presupposes a deeper look at the real basis for the results, at structure beyond size – which would not be bad.

It is not unstructured size that sets the real numbers apart from the rational numbers. The real numbers are different from the rational numbers for exactly the same reason that the rational numbers are different from the natural numbers: Not because there are 'more', but because their ordering is different. The real numbers thus remain distinguished from the rational numbers by their continuous ordering even when they are not assigned a higher cardinality.

When substantive results in mathematics – results outside the hedged garden of cardinal arithmetic – rely solely on size without any basis in ordering, topology, algebraic or other tangible structure then they must be suspected of having been manufactured on the cheap. To give an example: There are good reasons for believing in the existence of transcendental numbers. The argument from the countability of the algebraic numbers compared to the uncountability of all reals is not one of them.

In short, there are gains in both directions from carrying out the programme: Where traditional postulates that imply uncountability are maintained, they will at least have to be declared openly and justified against equally viable alternatives. Where single infinity is adopted, results that may have been assumed too glibly in the past will require renewed rewarding work, in a closer engagement with structure.

People are free to play with any formalism they choose, and in a liberal society it is no one's business to prohibit harmless activities that involve the postulation of multiple infinities. What does bear saying, however, is that the bent world postulates appear to be over-researched. And if indeed they are exhausted the restructuring programme which began with logic and the foundations needs to be extended far into mathematics proper.

# Second Order Logic

Inside a theory, proving theorems, there has never been an effective difference between the first order induction scheme

For all predicate strings P, P(0) $\wedge$ ( $\forall$x P(x) $\rightarrow$ P(s(x )) ) $\rightarrow$ $\forall$x P(x) is an axiom

and the second order induction axiom

$\forall$P P(0) $\wedge$ ( $\forall$x P(x) $\rightarrow$ P(s(x)) ) $\rightarrow$ $\forall$x P(x).

Scheme and axiom allow precisely the same derivations provided only that the same predicate strings are allowed to be substituted into either.

**Thesis:** The definition of order for a theory ought to be based on the range of substitution.

Because of the interchangeability of axiom and scheme, what matters is quantification in predicates, not axioms: First order number theory is the theory where only predicates with bound first order variables are allowed to be substituted into either scheme or axiom; in second order number theory the range of substitution extends to predicates containing bound second order variables.

There is a difference between axiom and scheme when mapping theories. The scheme will only achieve the same effect as $\forall$P in a model if the modelling transformation is onto. So although there is otherwise little to choose between them we propose to adopt the single axiom as the more honest way of doing what the infinite scheme was intended to do.

**Thesis:** Axiomatisability ought to be defined as finite axiomatisability.

Insisting on finite axiomatisability makes no difference for theories that can alternate between an infinite scheme and a single axiom. But it imposes necessary discipline on the ideas for theories that cannot so alternate. It helps avoid infinite lists and the unfortunate symbol '...'

The infinite list of sentences

0$\neq$d, 1$\neq$d, 2$\neq$d, ...

is apparently not reducible to a single axiom, or at least only to an axiom employing unconventional notation:

"$\forall$n s$^{(n)}$ $\neq$ d"

Second order *theories* are, and have always been, compact. Second order *logic* was thought to be non-compact because the models traditionally used for both first and second order logic have tolerated forms of inference that give rise to non-compactness.

Traditional models implicitly relied on the rule of inference: P(0), P(0'), P(0''), ... $\Rightarrow$ $\forall$x P(x).

Or rewritten without triple dots: $\forall$n P(s$^{(n)}$) $\Rightarrow$ $\forall$x P(x).

It is not at all clear how the sentence $\forall$n P(s$^{(n)}$) was meant to be derived for a model. Definite rules for this form of quantification inside terms were never articulated; the meaning always remained obscure. Most of the extant reasoning flees from numerals to vague talk of numbers, retreating into "extensional intuition".

No theory can actually support derivations from infinitely many premises. There is no agent to carry out such derivations. What a theory can do is to allow unconventional forms of notation and/or derivation that could be addressed, in a metaphorical way, as having effected an infinite derivation, when in fact the derivation only passes from one finitely long string to another. No matter where it appears infinity is only a simulation created with finitely many strings.

**Thesis:** Theories must be compact.

All logic is compact. What was discussed under the heading of "non-compactness" concerns the unusual properties theories acquire when supplemented by unconventional notation such as the quantifier in "$\forall$n P(s$^{(n)}$)".

The pyramid of sentences

$$\exists x \; x \neq 0$$

$$\exists x \; x \neq 0 \wedge x \neq 1$$

$$...$$

or equivalently its summary in unconventional notation

$$\forall n \; \exists x \bigwedge_{i<n} x \neq s^{(i)}$$

are "non-compact" in the sense that although every sentence in the pyramid scheme could be provable in a model, the limit

$$\text{"} \exists x \; \forall n \; x \neq s^{(n)} \text{"}$$

need not. (The limit only exists when quantifiers can be swapped to pass from $\forall n \exists x$ to $\exists x \forall n$ – clearly not a move that is valid in general.)

A similar diagonal argument can be made for the second order:

$$\exists P \; P(0) \leftrightarrow \neg P_0(0)$$

$$\exists P \; P(0) \leftrightarrow \neg P_0(0) \wedge P(1) \leftrightarrow \neg P_1(1)$$

$$...$$

There is finally no good reason to believe that second order quantification is somehow less sound than unrestricted first order quantification. The real hiatus is between finite quantification $\forall n < N$ and unrestricted quantification $\forall n$, not between $\forall n$ and $\forall P$ or $\forall f$, all three of which are all essentially and ineliminably infinitistic. If there is anything that is suspect and needs to be handled with care it is the extensional forms of reasoning that used to figure in traditional models. A rigorous treatment would replace sentence pyramids with the predicate sequences presented in the Cantor section.

## Non-Standard Models

For theories, the idea of an extension is simple – add more sentences. For models, i.e. theories under a transformation, a distinction needs to be drawn between two ways in which sentences can be added.

> **Definition:** A theory M is said to refine a theory Th if M is a model for Th under a non-surjective transformation.

> **Definition:** A theory E is said to enlarge a theory Th if a proper subtheory of E (but not E itself) is a model for Th.

The difference between refinement and enlargement is that of being a model as opposed to containing a model. Let $\mathbb{R}\wedge\mathbb{Q}$ be the theory of the ordering and algebraic properties shared by the rational and real numbers. $\mathbb{R}\wedge\mathbb{Q}$ can be refined in two directions, by an axiom saying that every element is a rational fraction to yield the structure of the rational numbers, or by a completeness axiom for sequences to yield the real numbers. $\mathbb{R}$ is then an enlargement of $\mathbb{Q}$, but not a refinement; $\mathbb{Q}$ is neither an enlargement nor a refinement for $\mathbb{R}$.

> **General Principle:** Every theory has non-isomorphic enlargements.

There is always more to add to the menagerie. A popular if unimaginative way of manufacturing enlargements is given by

> **Definition:** The Cartesian product of two theories $Th_1$ and $Th_2$ is a theory consisting of $Th_1$ relativised to a fresh predicate letter $U_1$ and $Th_2$ relativised to a fresh predicate letter $U_2$.

The general principle would hardly be worth mentioning had not a special case been held up as an important theorem: Let the theory Th be extensionalisable with an infinite universe u. Enlarge its model in set theory by the power set of its universe. The only constraint on the cardinality of the power set is then given by $|V| \geq |p(u)| \geq |u| \geq |\omega|$. In a standard shot, all four cardinalities will be equal. In different non-standard shots $|p(u)|$ could be set equal to any cardinality upwards of $|u|$.

> **The Löwenheim-Skolem Theorem:** Every extensionalisable theory with an infinite universe has enlargements that are non-isomorphic because of the size assigned to the universe and its power set.

A theory can constrain the elements of its universe, and the elements of the power set of its universe. The universe and the power set itself are not objects of the theory (Set theory is no exception to this rule since it only manages to talk about its own universe by becoming non-extensionalisable). The relation of the universe to its power set is therefore not a relation internal to any theory. In particular, the existence or otherwise of bijective functions between them must be left open. As the theory does not rule on it, almost anything is possible in enlargements.

Concerning set theory, the Löwenheim-Skolem Theorem applies only to subtheories with a suitably restricted universe.

> **Definition:** A set over which equality and membership are extensionalisable is called ordinary.

Initial segments of the constructible hierarchy, including most common mathematical sets, are ordinary. V is not.

To summarise up to this point: The Löwenheim-Skolem Theorem is not a source of what were called "non-standard" models because the models it says exist are enlargements. Non-standard models properly so called would have to be refinements.

If the underdetermination of number theory that would allow it to have non-isomorphic models was partly illusory, the tangible theory of non-standard arithmetic that was inspired by it was not. This theory is just a somewhat modified set of axioms not for number theory, but for – non-standard arithmetic.

> **Proposition:** A theory is categorical iff it has no refinements.
> Proof: Immediate.

First order number theory is not categorical. It has second order number theory and a range of non-standard theories as refinements.

Say that two numbers are in the same segment if they are finitely many successor steps away from each other. Standard number theory is then characterised by the fact that there is only one segment, non-standard arithmetics by any number of segments other than one. It of course takes either set theory or second-order predicates to define segments. (The standard refinement is completed simply by allowing the previously unavailable predicate P(n) = "n is finitely many successor steps away from 0" to feed into induction; in non-standard theories induction must be modified so that P(n) does not become true of all numbers.)

The general way forward: Instead of questionably "constructing" non-standard models, define a theory by modified axioms. If these modified theories should turn out to be modelled by set theory, then that is a bonus, but never a requirement.

# Measure Theory

Measure Theory as we currently know it is valid only in the bent world; measure theory under single infinity has to be done differently.

**Proposition:** Under single infinity, no measure can be infinitely additive for general sets.
Proof: Infinite additivity would equate the measure of the continuum with that of an isolated point.

What might be conceivable is that a measure remains infinitely additive for certain well-behaved, non-zero measure sets. I state it as an open issue (open at least to me) how far the meaning of 'well-behaved' can be expanded beyond the obvious – beyond finite unions and intersections of connected sets.

In the bent world, sets of measure zero were defined to include all *countable* sets. That cannot continue in the standard world.

**Definition (untested):** A set is said to be of measure zero if it vanishes after finitely many derivations.

It may be possible to reconstruct Lesbesgue integration based on this changed definition of a negligible set. Lesbesgue's still brilliant idea of defining integration not for individual functions but rather for functions modulo zero measure sets, that is equivalence classes of functions that differ in value over at most a set of measure zero, could be maintained but with a different meaning given to 'zero measure'.

**Corollary:** No dense sets are measure zero.

In particular, the indicator function of the irrationals $\mathbb{R}\backslash\mathbb{Q}$ over the unit interval [0;1] no longer counts as integrable.

**Illustration:** Given some $\varepsilon > 0$, a dense set s, and instructions to draw a super fine vertical line of unit height and width less than $\varepsilon$ for every point in the set, and no lines for points not in the set, then which picture to draw, the one on the left or the one on the right, would be undecidable on the basis of the instructions.

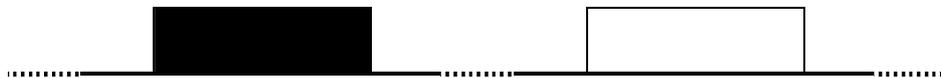

So that integrals are no longer defined for functions discontinuous over a dense set is perhaps as it should be. I have always found it intuitively offensive to assign the value 1 (all black) to the indicator integral over the dense set of irrationals $\mathbb{R}\backslash\mathbb{Q}$, but the value 0 (all white) to the indicator integral of the equally dense set of rationals $\mathbb{Q}$. Any value in between (shade of grey) would have been just as justified, and just as arbitrary.

# Topology

A major source of non-standard assumptions in topology are the postulates underwriting infinitary constructions. It may seem ironic that constructions, which were advocated by intuitionists and other reformers as particularly "safe", should be a cause for concern. But this is not surprising once it is considered how the nature of mathematical construction has been misunderstood.

The trouble with constructive arguments is that they are often nothing of the sort. That existence is being posited, not proven. And the word 'construction' allowed to imply success, when success is at issue. Although not as beguilingly intuitive, a non-constructive argument at least always constructs something real: a proof by contradiction, which is a well-defined finite string object. A constructive argument on the other hand might fundamentally misconstrue the business of construction – by, for instance, attempting object constructions that are not grounded in constructions with strings –, and conjure up a miasma.

Constructions do not, by themselves, prove existence; constructions are applications of postulates of existence. The definition of the object under construction is inevitably a string produced according to certain rules. Only when combined with a postulate – better declared openly but now often implicit – that all strings defined in accordance with this method are non-empty will the "construction" have succeeded in producing an object. Without the postulate there would only be a name, and no reason not to believe it to be empty.

Cantor's uncountability argument is the classic case of a construction that depends on rarely acknowledged postulates. A similar if less famous construction attributed to the same author, with the same ultimate consequence of uncountability, is the Cantor Intersection Theorem. (The Cantor Intersection Theorem, also known as the Chinese Box Theorem, serves as a lemma in the proof of the Heine-Borel Theorem, which then implies the uncountability of the real numbers.) The success of both constructions becomes refutable under single infinity.

For general compact spaces the Cantor Intersection Theorem is false. The construction used can push out only finitely far – unless the following assumption is added: That the initial segments, whose constructability may indeed be granted, are more than just initial segments for a sequence.

Initial segments for a sequence exist, initial segments for a sequence that would converge to a sequence of initial segments exist, etc. But assembling these finite segments into a whole requires a postulate.

The Cantor Intersection Property – that certain unusual sequences of nested closed sets exist and their intersection is non-empty – therefore fails to be provable outright. Where the property is desired, lay either it down as an axiom, or a construction principle from which it follows. Be aware, however, that both postulates lead to non-standard spaces.

To rescue the Heine-Borel Theorem for standard spaces its premise may have to be strengthened from

*If a sequence S of open intervals covers [0;1]* then some finite subset of S covers [0;1].

To

*If S is a sequence of open intervals such that it covers [0;1] and the midpoints of the intervals vanish after finitely many derivations* then ...

It is as yet unclear if reconstruction along these lines will prove viable.

## Odds and Ends

To close, an additional number of theorems that I suspect ripe for restructuring. The list does not claim to be complete, or even correct.

- From (Meta-)Set Theory: Forcing
  (Uses Gödel-numbering for infinite string sets, hence implicitly concatenation)

- From Topology: The Baire Category Theorem
  (A consequence of the Cantor Intersection Theorem.)

- From Cardinal Arithmetic: $|\omega \times \omega| = |\omega|$.
  (Cantor's image – an enumeration winding its way antidiagonally across a table – suggests much, but proves little. It does not preclude the possibility of defining a version of set theory in which the cardinality of a Cartensian product is larger than the cardinality of its factors.
  The only 'countable' sets, in any real sense relating to counting and ability, are finite. The remaining technical meaning of 'countable' refers to the existence of a bijective function beween the set and $\omega$. These functions are infinitary objects, which only a postulate can introduce. No construction, without relying on a postulate, could vouch for their existence.)

- From Model Theory: "The (countable) categoricity of a dense linear order."
  (Another theorem of Cantor's. What the construction establishes is probably only that any finite linearly ordered set can be mapped order-isomorphically into every interval of a dense set. The statement is equivalent to denseness, and could serve as another definition of dense.)

- From Model Theory: Łos' Theorem

## Afterword

No single paper could ever complete a restructuring programme, nor could a programme be completed by a single author. Restructuring, by its nature, requires a collective effort by many active in a discipline. I have suggested an outline, and hope that others will be able to see something attractive in that outline. I do not claim that everything contained in this manuscript is absolutely right or fully on course. I only claim that it is broadly pointing in a more promising direction than the status quo of the discipline. What follows, if indeed anything does follow, is up to others.

## "References"

Since no citations were made, there are also no references, just acknowledgements. Whereas credit for the ideas that were only restructured continues to belong to their first authors (not listed here), the original ideas are indebted especially to the following publications: